\documentclass{article}
\usepackage{amsfonts}
\begin{document}
\newtheorem{Pretheorem}{Pretheorem}{\bf}{\it}
\def\R{{\mathbb R}}
\def\C{{\mathbb C}}
\def\Z{{\mathbb Z}}
\def\H{{\mathbb H}}
\newcommand{\Tr}{\mbox{\rm Tr}}
\renewcommand{\Re}{\mbox{\rm Re}}
\renewcommand{\Im}{\mbox{\rm Im}}
\newtheorem{lemma}{Lemma}
\newtheorem{theorem}{Theorem}
\newtheorem{definition}{Definition}
\newtheorem{corollary}{Corollary}
\newtheorem{proposition}{Proposition}
\newtheorem{statement}{Statement}

\title{Dirac operators and conformal invariants of tori in
$3$-space}
\author{Iskander A. TAIMANOV
\thanks{Institute of Mathematics,
    630090 Novosibirsk, Russia,
    e-mail: taimanov@math.nsc.ru}}
\date{}

\maketitle

\section{Introduction}

In this paper we show how to assign to any torus immersed into
the three-space $\R^3$ or the unit three-sphere $S^3$ a complex
curve such that the immersion is described by functions defined on this
curve (a Riemann surface which is generically of infinite genus).
We call this curve the spectrum of a torus (with a fixed conformal
parameter).
This spectrum has many interesting properties and, in particular,
relates to the Willmore functional whose value is encoded in it.

The construction of the spectra for tori in $\R^3$ was briefly explained
in \cite{T3}. In this text we do that also for immersed tori in $S^3$.

Our conjecture that the spectrum of a torus in $\R^3$
is invariant under conformal transformations of $\R^3$
was proven modulo some analytic facts by Grinevich
and Schmidt \cite{GS}.
In fact, their proof is rather physical which may be expected because the
construction of the spectrum originates in soliton theory.

In this text we give a complete proof of the
conformal invariance of the spectra for isothermic tori.
This case already covers many interesting surfaces such as constant mean
curvature tori and tori of revolution in $\R^3$.

Some spectral curves of finite genus already appeared
in studies of harmonic tori in $S^3$ by Hitchin \cite{Hitchin} and
constant mean curvature (CMC) tori in $\R^3$ by Pinkall and Sterling
\cite{PS}. It was shown that such tori are expressed in terms of
algebraic functions corresponding to these complex curves
\cite{Hitchin,Bo1}. We show that for minimal tori in $S^3$
and CMC tori in $\R^3$ these spectral curves
are particular cases of the general spectrum.

The general construction is based on the global Weierstrass
representation of closed surfaces introduced in \cite{T1,T3}
and a general construction of the Floquet(or Bloch) variety for a periodic
differential operator. An existence of this variety is derived from
the Keldysh theorem but an effective construction which gives more
information about analytic behavior of this complex curve was proposed by
Krichever in \cite{Kr} who used perturbation methods.
It is as follows. Take an immersed torus with the induced metric
$e^{2\alpha} dz d\bar{z}$ and consider differential operators
$$
{\cal D} =
\left(\begin{array}{cc}
0 & \partial \\ -\bar{\partial} & 0
\end{array} \right) +
\left(\begin{array}{cc}
U & 0 \\ 0 & U
\end{array} \right)
\ \ \ \ \mbox{with $U=\frac{1}{2}He^{\alpha}$ for a torus in $\R^3$},
$$
$$
{\cal D}^S =
\left(\begin{array}{cc}
0 & \partial \\ -\bar{\partial} & 0
\end{array} \right) +
\left(\begin{array}{cc}
V & 0 \\ 0 & \bar{V}
\end{array} \right)
\ \ \ \ \mbox{with $V=\frac{1}{2}(H-i)e^{\alpha}$ for a torus in $S^3$},
$$
where $H$ is the mean curvature.

Let $\Lambda$ be the period lattice of a torus which means that the torus
is an immersion of $\C/\Lambda$ with a conformal parameter $z \in \C$ on it.
Take a basis $\gamma_1,\gamma_2$ for
$\Lambda$ which is also considered as a basis for $H_1(T^2) \approx \Lambda$.
Now consider all solutions $\psi$ to the equations
$$
{\cal D}\psi = 0 \ \ \ \ \ \mbox{or} \ \ \ \ \ {\cal D}^S\psi = 0
$$
satisfying the following conditions
$$
\psi(z + \gamma_j) = \mu_j \psi(z).
$$
These are Floquet(--Bloch) functions and the pairs $(\mu_1,\mu_2)$
form a complex curve in $\C^2$. This is the Floquet
zero-level spectrum of ${\cal D}$ and, by the definition,
this is the spectrum of the immersed torus.
The analytic properties of this curve are described by Pretheorem
which is a modification of its analogs for two-dimensional scalar Schr\"odinger
and heat operators proven in \cite{Kr} and it is clear that the proof of
Pretheorem may be obtained by slight modifications of
the reasonings of \cite{Kr}.

One of the most interesting properties of this construction is its relation
to a conformal geometry and the Willmore functional which
equals
$$
4 \int_{\C/\Lambda} U^2 dx \wedge dy \ \ \ \mbox{or} \ \ \
4 \int_{\C/\Lambda} |V|^2 dx \wedge dy
$$
for tori in $\R^3$ or $S^3$, where $z = x+iy$.

The global Weierstrass representation of closed
surfaces represents any
closed surface $\Sigma$ in terms of a solution
to the equation ${\cal D}\psi = 0$ (a harmonic spinor) where
$\psi$ takes values in some bundle over
the constant curvature surface $\Sigma_0$ which is
is conformally equivalent to $\Sigma$ (see Theorems 1--3 in 2.3)
\cite{T1,T3}.
The Willmore functional $\int_{\Sigma} (H^2 - K)d\mu =
4 \int_{\Sigma_0} U^2 d x \wedge d y - 2\pi \chi(\Sigma)$
measures the $L_2$-norm of the potential $U$ of the surface.
For small values of this functional the equation
${\cal D}\psi = 0$ does not admit solutions which describe closed surfaces in
$\R^3$ and therefore that explains physical meaning
of lower bounds for the Willmore
functional proposed by the Willmore conjecture and its generalizations.
This gives a hint that the spectral properties of ${\cal D}$ have to have
a geometric meaning. We shall discuss that in details in 4.4.

In \cite{T4} it was established that the dimension of the kernel of ${\cal D}$
gives lower estimates for the Willmore functional for spheres.
Actually we proved for spheres with one-dimensional
potentials $U$ (examples of them are spheres of revolution but not only)
and conjectured for all spheres the following inequality
\begin{equation}
\int_{\Sigma_0} U^2 d x \wedge d y \geq 4 \pi N^2
\ \ \ \ \mbox{with $N =  \dim_{\C} \ker {\cal D}$.}
\label{wn}
\end{equation}
The proof of it
is based on the inverse scattering problem for the one-dimensional
Dirac operator and it also works for general Dirac operators
with $S^1$-symmetry (i.e., with one-dimensional potentials)
on special spinor bundles over the 2-sphere. This inequality can not be
improved and the equality is achieved on ``soliton spheres'' \cite{T4}.

Another treatment of the global representation belongs
to Pedit and Pinkall who proposed to consider spinor $\C^2$-bundles
introduced in \cite{T1} as quaternionic line bundles
and consider harmonic spinors as holomorphic quaternionic sections of
such bundles. This enables them to apply ideas of algebraic
geometry to surface theory and to generalize this representation for
surfaces in $\R^4$ \cite{PP}. Very recently
they managed to relate (\ref{wn}) to the quaternionic analog of the
Pl\"ucker formula and by that prove our conjecture, i.e.,
establish  the inequality (\ref{wn}) for all spheres
together with its generalizations for higher genus surfaces.

This paper is organized as follows.

In section 2 we recall the notion of
the Weierstrass representation.

In section 3 we prove that the multipliers
$(\mu_1,\mu_2)$ of Floquet functions form a spectral curve in
$\C^2$ and discuss its analytic properties.

In section 4 we show how
to assign such a spectrum to an immersed torus in $\R^3$, show the Willmore
functional appear in this picture and discuss the modern state of the Willmore
conjecture.

In section 5 we show that the spectra of CMC and isothermic tori
are particular cases of the spectrum defined in section 4 and also prove our
conjecture that the spectra of an isothermic torus and its dual surface
coincide \cite{T2}.

In section 6 we show how to
assign a spectral curve to a torus in $S^3$ and prove that for
a minimal torus in $S^3$ it coincides with
a spectral curve defined by Hitchin for harmonic tori in $S^3$
\cite{Hitchin}.

In section 7 we prove that
the spectrum of an isothermic torus in $S^3$ is invariant
with respect to conformal transformations of $\bar{\R}^3 =
\R^3 \cup \{\infty\}$
preserving the torus in $\R^3$.

\section{The Weierstrass representation}

\subsection{Basic equations of surface theory in $\R^3$}

In this subsection we recall the main definitions and some well-known
important facts from the classical surface theory.

Let ${\cal U}$ be a domain in $\R^2$, with coordinates $(x^1,x^2)$,
regularly immersed into $\R^3$:
$$
F: {\cal U} \rightarrow \R^3.
$$
At every point $p \in {\cal U}$ the vectors
$$
F_1 = \frac{\partial F}{\partial x^1}, \ \ \
F_2 = \frac{\partial F}{\partial x^2}, \ \ \
N=\frac{[F_1 \times F_2]}{|F_1||F_2|}
$$
form a linear basis $\sigma = (F_1,F_2,N)^{\top}$
for $\R^3$, where $F_1$ and $F_2$ are tangent vectors to the surface $\Sigma =
F({\cal U})$,
and $N$ is a unit normal vector.
The variables $(x^1,x^2)$ are local coordinates on
$\Sigma$ and the induced metric on it is
$$
{\bf I} = g_{kl}dx^k dx^l \ \ \mbox{with $g_{kl} =
\langle F_k,F_l \rangle$ \ \ (the {\it first fundamental form})}
$$
where $\langle a, b \rangle = a_1 b_1 + a_2 b_2 + a_3 b_3$.

The derivatives of the basic vectors are expanded in $\sigma$ as
\begin{equation}
\frac{\partial^2 F}{\partial x^k \partial x^l} =
\Gamma^j_{kl}\frac{\partial F}{\partial x^j} +
b_{kl}N, \ \ \
\frac{\partial N}{\partial x^k} = - b^j_k\, \frac{\partial F}{\partial x^j}
\label{gw}
\end{equation}
where $\Gamma^j_{kl}$ are the {\it Christoffel symbols},
${\bf II} = b_{kl} dx^k dx^l$ is the {\it second
fundamental form}, and $b^j_k = g^{j l}b_{l k}$.
The equations (\ref{gw}) are the {\it Gauss--Weingarten derivation
equations} and have the form
\begin{equation}
\frac{\partial \sigma}{\partial x^1} = {\bf U}\sigma, \ \ \
\frac{\partial \sigma}{\partial x^2} = {\bf V}\sigma,
\label{derivation}
\end{equation}
where $U$ and $V$ are $(3\times 3)$-matrices. The compatibility conditions for
(\ref{derivation}) are the {\it Codazzi equations}:
$$
\frac{\partial^2 \sigma}{\partial x^1 \partial x^2} -
\frac{\partial^2 \sigma}{\partial x^2 \partial x^1}
=
\left(\frac{\partial {\bf U}}{\partial x^2} -
\frac{\partial {\bf V}}{\partial x^1} +
[{\bf U},{\bf V}]\right) \sigma = 0,
$$
which are equivalent to the {\it zero-curvature equations}
\begin{equation}
\frac{\partial {\bf U}}{\partial x^2} -
\frac{\partial {\bf V}}{\partial x^1} +
[{\bf U},{\bf V}]= 0
\label{codazzi}
\end{equation}
for the connection $(\partial/\partial x^1 -
{\bf U},\partial/\partial x^2 - {\bf V})$.

At every point $p$ of the surface the fundamental forms are diagonalized as
$$
{\bf I}(p) = \left(
\begin{array}{cc} 1 & 0 \\ 0 & 1
\end{array}
\right),\ \ \
{\bf II}(p) = \left(
\begin{array}{cc} k_1 & 0 \\ 0 & k_2
\end{array}
\right)
$$
and the {\it principal curvatures}
$k_1$ and $k_2$ satisfy the equation
$\det\left(b_{kl} - k g_{kl}\right) = 0$,
which divided by $\det g_{jk}$ takes the form
$k^2 - 2H k + K = 0$
where $H$ is the {\it mean curvature} and
$K$ is the {\it Gaussian curvature}:
$$
H = \frac{k_1 + k_2}{2}, \ \ \
K = k_1 k_2.
$$
A point $p$ is called an {\it umbilic} point if the principal
curvatures coincide at $p$: $k_1 = k_2$, which is equivalent to
$H^2 - K = 0$.

Let $z=x^1 +ix^2$ be a conformal parameter on the surface, i. e.,
the first fundamental form is
$$
{\bf I} = e^{2\alpha (z,\bar{z})} d z\,d\bar{z},
$$
which means that
$$
\langle F_z, F_z \rangle = \langle F_{\bar{z}}, F_{\bar{z}} \rangle = 0,
\ \ \ \langle F_z, F_{\bar{z}} \rangle = \frac{1}{2}e^{2\alpha}.
$$
The family
$\widetilde{\sigma} = (F_z, F_{\bar{z}}, N)^\top$ is a basis for $\C^3$ and
the Gauss--Weingarten equations are written as
$$
\widetilde{\sigma}_z = \widetilde{\bf U}\sigma, \ \ \
\widetilde{\sigma}_{\bar{z}}= \widetilde{\bf V}\sigma,
$$
with
$$
\widetilde{\bf U} =
\left(
\begin{array}{ccc}
2\alpha_z & 0 & A \\
0 & 0 & B \\
-2e^{-2\alpha}B & -2e^{-2\alpha}A & 0
\end{array}
\right),\ \
\widetilde{\bf V} =
\left(
\begin{array}{ccc}
0 & 0 & B \\
0 & 2\alpha_{\bar{z}}  & \bar{A} \\
-2e^{-2\alpha} \bar{A} & - 2e^{-2\alpha}B & 0
\end{array}
\right),
$$
$A = \langle F_{z z}, N \rangle$, and
$B = \langle F_{z \bar{z}}, N \rangle$.
The second fundamental form equals
$$
{\bf II} = (2B + (A + \bar{A}))(d x^1)^2 + 2i(A-\bar{A})d x^1 d x^2 +
(2B- (A + \bar{A})) (d x^2)^2,
$$
and we have
$$
H = 2Be^{-2\alpha}, \ \ \ K = 4(B^2 - A\bar{A})e^{-4\alpha}.
$$
Now the Codazzi equations
$\widetilde{\bf U}_{\bar{z}} - \widetilde{\bf V}_z +
[\widetilde{\bf U},\widetilde{\bf V}] = 0$
take the form
\begin{equation}
\alpha_{z\bar{z}} + e^{-2\alpha}(B^2 - A\bar{A}) = 0, \ \ \
A_{\bar{z}} - B_z + 2\alpha_z B = 0.
\label{conf-codazzi}
\end{equation}
The first of them is the {\it Gauss egregium theorem}
and another equation
$$
A_{\bar{z}} = \frac{1}{2} H_z\, e^{2\alpha}
$$
splits into two real-valued equations.

A quadratic differential $\omega = A\, d z^2$ is called
the {\it Hopf differential} and has  important geometrical
properties. For instance, $\omega$ vanishes at a point
if and only if this is an umbilic point.

It is said that the {\it Gauss map} of a surface $\Sigma = F({\cal U})$
$$
G: \Sigma \rightarrow S^2, \ \ \ \ G(p)=N(p),
$$
is {\it harmonic} if
$\Delta G (p) = \lambda(p) N(p)$, where
$\Delta$ is the {\it Laplace--Beltrami operator}:
$\Delta = 4 e^{-2\alpha}\partial \bar{\partial}$.
We have  $\Delta F = 2HN$.
Since $N_{z\bar{z}} = -2e^{-2\alpha}(\bar{A}_z F_z + A_{\bar{z}}F_{\bar{z}}
+(A\bar{A} + B^2)N)$, we conclude that

{\it
the Gauss map $G$ is harmonic if and only if the Hopf differential
$\omega$ is holomorphic, which, by (\ref{conf-codazzi}), is equivalent to
$H_z = H_{\bar{z}} = 0$, i.e., $H = \mbox{\rm const}$ and
$\Sigma$ is a constant mean curvature (CMC) surface} \cite{RV}.

There are two other important classes of surfaces:

1) a surface is called {\it minimal} if $H = 0$, which is equivalent to
$F_{z\bar{z}} = 0$;

2) a surface is called {\it isothermic} if there is a conformal parameter on
it such that $\mbox{\rm Im}\, A = 0$.

\subsection{The local representation of a surface}

In this subsection we follow \cite{Kon,T1}.

Denote by ${\cal Q}$ a quadric in $\C^3$ defined by the equation
$$
Z_1^2 + Z_2^2 + Z_3^2 = 0, \ \ \ \ Z=(Z_1,Z_2,Z_3) \in \C^3.
$$
For a conformal parameter $z$ on a surface $\Sigma$
there is a mapping
\begin{equation}
f:{\cal U} \stackrel{F}{\to} \Sigma \to {\cal Q} \ \ \mbox{where
\ $f(p) = F_z(p)$},
\label{map}
\end{equation}
satisfying the conditions
\begin{equation}
\Im \frac{\partial f}{\partial \bar{z}} = 0.
\label{new1}
\end{equation}
It is clear that
any mapping $f:{\cal U} \rightarrow {\cal Q}$ satisfying (\ref{new1})
has the form $f = \partial_z \Phi$ for some real-valued function
$\Phi$ and therefore has the form (\ref{map}) for some surface.

The set ${\cal Q}$ is parameterized by
$(\varphi_1,\varphi_2) \in \C^2$ as follows:
$$
Z_1=\varphi_1^2 - \varphi_2^2,
\ \ Z_2=i(\varphi_1^2 + \varphi_2^2),
\ \ Z_3 = 2\varphi_1\varphi_2.
$$
For simplicity, renormalize $\varphi$ as follows:
$$
\psi_1 = \sqrt{\frac{2}{i}}\varphi_1, \ \
\psi_2=\sqrt{\frac{2}{i}}\bar{\varphi}_2.
$$
In terms of $\psi$ the equations (\ref{new1}) are
\begin{equation}
{\cal D} \psi = 0
\label{diracequation}
\end{equation}
where
\begin{equation}
{\cal D} =
\left( \begin{array}{cc} 0 & \partial \\ -
\bar{\partial} & 0 \end{array} \right) +
\left(
\begin{array}{cc}
U & 0 \\
0 & U
\end{array}
\right)
\label{diracoperator}
\end{equation}
is a Dirac operator with a real-valued potential $U(z,\bar{z})$.

Now, if $F(p_0) = (x^1_0,x^2_0,x^3_0) \in \R^3$, then the surface is
described by the {\it Weierstrass formulas}
$$
x^1(p) = x^1(p_0) + \int_{p_0}^p \left(
\frac{i}{2} (\bar{\psi}_2^2 + \psi_1^2) d z -
\frac{i}{2} (\psi_2^2 + \bar{\psi}_1^2) d\bar{z}\right),
$$
\begin{equation}
x^2(p) = x^2(p_0) + \int_{p_0}^p \left(
\frac{1}{2}(\bar{\psi}_2^2 - \psi_1^2) d z +
\frac{1}{2}(\psi_2^2 - \bar{\psi}_1^2) d\bar{z} \right),
\label{weierstrass}
\end{equation}
$$
x^3(p) = x^3(p_0) + \int_{p_0}^p \left(
\psi_1 \bar{\psi}_2 d z + \bar{\psi}_1 \psi_2 \right),
$$
which are just
$$
F(p) = F(p_0) + \int_{p_0}^p (f d z + \bar{f} d \bar{z}).
$$
These local formulas
in different forms were known before (see, for instance,
\cite{Ken,Abresch,Bo1} and comments in \cite{T1,T3}) but this form was
introduced by Konopelchenko \cite{Kon} who working with the formulas of
Eisenhart \cite{Eis} elaborated them into a form most convenient for
applications. He considered them for constructing some surfaces via
solutions to (\ref{diracequation}) and defining soliton deformations of
``induced surfaces'' but as we see this gives a general
local construction of surfaces.

By straightforward computations it is derived that
\begin{equation}
U = \frac{H e^{\alpha}}{2},
\ \ e^{\alpha} = |\psi_1|^2 + |\psi_2|^2.
\label{formulas}
\end{equation}
and we see that for $H=0$ the formulas (\ref{weierstrass}) reduce to
the classical formulas for minimal surfaces.

Compute
$N = e^{-\alpha} \left(
-i(\bar{\psi}_1 \bar{\psi}_2 - \psi_1 \psi_2),
-(\bar{\psi}_1 \bar{\psi}_2 + \psi_1 \psi_2), |\psi_2|^2-|\psi_1|^2
\right)$
and taking (\ref{diracequation}) into account derive
\begin{equation}
A = \langle F_{z z}, N \rangle = \psi_{1 z}\bar{\psi}_2 -
\bar{\psi}_{2 z}\psi_1, \ \
B = \langle F_{z\bar{z}}, N \rangle = U e^{\alpha}.
\label{A}
\end{equation}

The Gauss--Weingarten equations written in terms of $\psi$
describe the deformations of $\psi$
and the first half of them is just the Dirac equation (\ref{diracequation}).
For obtaining another half of equations differentiate $e^{\alpha}$ by $z$:
$$
\alpha_z e^{\alpha} = \bar{\psi}_1\psi_{1 z} + \psi_2\bar{\psi}_{2 z},
$$
and, taking (\ref{A}) into account, obtain
$$
\psi_{1 z} = \alpha_z \psi_1 + Ae^{-\alpha}\psi_2,
\ \ \
\psi_{2\bar{z}} = -\bar{A}e^{-\alpha}\psi_1 + \alpha_{\bar{z}}\psi_2.
$$
Now the Gauss--Weingarten equations are written as
\begin{equation}
\left[\frac{\partial}{\partial z} -
\left(\begin{array}{cc}
 \alpha_z & A e^{-\alpha} \\
-U & 0
\end{array}
\right)\right]\psi =
\left[\frac{\partial}{\partial \bar{z}} -
\left(\begin{array}{cc}
0 & U \\
-\bar{A}e^{-\alpha} & \alpha_{\bar{z}}
\end{array}
\right)\right]\psi = 0
\label{w-gw}
\end{equation}
and the compatibility conditions for them,
the Codazzi equations, are
\begin{equation}
A_{\bar{z}} = (U_z - \alpha_z U)e^{\alpha}, \ \ \
\alpha_{z\bar{z}} + U^2 - A\bar{A}e^{-2\alpha} = 0.
\label{w-codazzi}
\end{equation}
In fact, the equation (\ref{diracequation}) is already the
compatibility condition
for an existence of a surface with the Gauss map given by $\psi$.
The other half of the equations (\ref{w-gw}) follows from it.

In a paper by Friedrich \cite{F}
this representation
was explained by classical means of the theory of Dirac operators.
We also would like to mention
a paper by Matsutani \cite{Mats} where it was
considered from the physical point of view.

\subsection{The global Weierstrass representation}

Here we follow \cite{T1} where this global representation was
introduced.

Any closed oriented surface $\Sigma$ in $\R^3$ is conformally equivalent to
a constant curvature surface $\Sigma_0$ and a choice of a conformal
parameter $z$ on $\Sigma$ means that a conformal equivalence
$\Sigma_0 \to \Sigma$ is fixed.

To define a compact oriented surface globally via the formulas
(\ref{weierstrass}) we have to introduce fibre bundles over surfaces
and Dirac operators on them.

Consider two cases:

1) {\sl Tori.}
Let $\Sigma$ be a torus immersed into $\R^3$. Then it is conformally
equivalent to a flat torus $\Sigma_0 = \C/\Lambda$ and $z$ is a conformal
parameter. The vector function $\psi$ is expanded to a section of
a $\C^2$-fiber bundle over $\Sigma$ defined by the monodromy rules
\begin{equation}
\psi(z + \gamma) = \varepsilon(\gamma)\psi(z)
\label{psi-torus}
\end{equation}
where $\gamma \in \Lambda$ and $\varepsilon: \Lambda \to \{\pm 1\}$ is a
character of $\Lambda$, i. e., a homomorphism to $\{\pm 1\}$.
The Dirac operator ${\cal D}$ acts on this bundle and
\begin{equation}
U(z + \gamma) = U(z).
\label{U-torus}
\end{equation}
Hence, we have

\begin{theorem}
$($\mbox{\rm{\cite{T1}}}$)$ \
The formulas (\ref{psi-torus}) and (\ref{U-torus}) define
a $\C^2$ bundle over a flat torus $\Sigma_0$. To any section $\psi$,
of this bundle, satisfying the Dirac equation (\ref{diracequation})
corresponds a surface in $\R^3$ defined by the formulas (\ref{weierstrass})
up to translations in $\R^3$.
\end{theorem}

2) {\sl Surfaces of genus $g > 1$.}
Let $\Sigma_0$ be a hyperbolic surface conformally equivalent
to a surface $\Sigma$ immersed into $\R^3$ and $z$ be a conformal parameter.
The surface $\Sigma_0$ is isometric to
${\cal H}/\Lambda$, where ${\cal H}$ is the Lobachevsky
upper-half plane and $\Lambda$ is a discrete subgroup of $PSL(2,\R)$.

Any element $\gamma \in \Lambda \subset
PSL(2,\R)$ is represented by elements
$$
\pm
\left(
\begin{array}{cc}
a & b \\ c & d
\end{array}
\right), \ \ \ \
a,b,c,d \in \R, \ \ \ ad-bc =1.
$$
The action on ${\cal H}$ is
$$
z \to \gamma(z) = \frac{az + b}{cz + d}.
$$

Define over $\Sigma_0$ a $\C^2$-bundle by the monodromy rules
\begin{equation}
\psi_1(\gamma(z)) = (cz+d)\psi_1(z), \ \ \
\psi_2(\gamma(z)) = (c \bar{z} + d)\psi_2(z).
\label{psi-handle}
\end{equation}
The Dirac operator acts on this bundle and
\begin{equation}
U(\gamma(z)) = |cz+d|^2 U(z).
\label{U-handle}
\end{equation}

\begin{theorem}
$($\mbox{\rm{\cite{T1}}}$)$ \
The formulas (\ref{psi-handle}) and (\ref{U-handle}) define
a $\C^2$-bundle over a hyperbolic surface $\Sigma_0$.
To any section $\psi$,
of this bundle, satisfying the Dirac equation (\ref{diracequation})
corresponds a surface in $\R^3$ defined by the formulas (\ref{weierstrass})
up to translations in $\R^3$.
\end{theorem}

Notice that $\psi_1 \sqrt{d z}$ and $\bar{\psi}_2 \sqrt{d z}$ are defined
modulo $\pm 1$. In fact, they are spinors and
therefore we shall call $\psi_1$ and $\psi_2$ also spinors.
There are $2^{2g}$ such nonequivalent spinor bundles over $\Sigma_0$
where $g$ is the genus of $\Sigma_0$.

This representation of compact oriented surfaces in $\R^3$ via
solutions of Dirac equations (i.e., harmonic spinors) in spinor bundles
over constant curvature surfaces is called the {\it Weierstrass
representation} of surfaces.

\begin{theorem}
$($\mbox{\rm{\cite{T1}}} for real-analytic surfaces,
\mbox{\rm{\cite{T3}}} for $C^3$-regular surfaces$)$
Every smooth closed oriented surface in $\R^3$ has a Weierstrass
representation.
\end{theorem}

We see from the direct construction in 2.2
that for a surface with a fixed conformal parameter such a
representation is unique. This gives rise to the following definition.

\begin{definition}
Let $(\Sigma,z)$ be an immersed surface with a fixed conformal parameter.
Then the potential $U$ of its Weierstrass representation is called
the potential of a surface.
\end{definition}

To any harmonic spinor on $\Sigma_0$ with a potential $U$
there corresponds a surface whose Gauss map descends
through $\Sigma_0$.
The criterion of closedness of such a surface is as follows.

\begin{proposition}
A surface represented by a harmonic spinor $\psi$ over a compact surface
$\Sigma_0$ is closed if and only if
$$
\int_{\Sigma_0} \bar{\psi}_1^2 \,
d\bar{z} \wedge \omega = \int_{\Sigma_0} \psi_2^2 \, d\bar{z}
\wedge \omega = \int_{\Sigma_0} \bar{\psi}_1 \psi_2 \, d\bar{z}
\wedge \omega = 0
$$
for any holomorphic differential on $\Sigma_0$.
\end{proposition}

This proposition was proved by M. Schmidt for general tori (in this case
$\omega = {\rm const}\cdot dz$) and by the author \cite{T4}
for higher genera surfaces.

One of the most important properties of this representation is the
equality
$$
4 \int_{\Sigma_0} U^2 d x \wedge d y = \int_{\Sigma} H^2 d \mu
$$
where $d \mu$ is the measure given by the induced metric on $\Sigma$
\cite{T1}. The functional
$$
{\cal W}(\Sigma) =
\int_{\Sigma} (H^2 - K) d \mu = \int_{\Sigma}\left( \frac{k_1-k_2}{2}
\right)^2 d\mu
$$
is called the {\it Willmore functional}.
By the Gauss--Bonnet theorem, for closed oriented surfaces
it equals
$$
{\cal W}(\Sigma) = \int_{\Sigma}H^2 d \mu - 2\pi \chi(\Sigma)
$$
where $\chi(\Sigma)$ is the Euler characteristic of the surface and,
therefore, for tori ${\cal W} = \int_{\Sigma}H^2 d \mu$.

There is a famous Willmore conjecture that for tori the Willmore functional
is greater or equal than $2\pi^2$. We shall discuss it in 4.4.

We see the main
advantage of the global Weierstrass representation in
using the spectral properties of ${\cal D}$
for study of conformal geometry of surfaces. The present paper is devoted to
developing this idea
for tori.

\section{The Floquet spectrum of a periodic operator}

\subsection{Floquet functions and the spectral curve}

Let $L$ be a differential operator acting on functions or vector functions
on $\R^n$, whose coefficients are
periodic with respect to translations by vectors from a lattice
$\Lambda$, i. e., $\Lambda$-periodic. We assume that $\Lambda$ has the maximal
rank, which means that it is isomorphic to $\Z^n$ and $\R^n/\Z^n$ is a torus.

To any vector $\gamma \in \Lambda$ there corresponds a translation operator
$T_{\gamma}$:
$$
T_{\gamma} f(x) \to f(x+\gamma),
$$
Since $L$ is $\Lambda$-periodic, if
$L f = \lambda f$, then $L T_{\gamma} f = \lambda T_{\gamma}f$.
Moreover
$$
[T_{\gamma},L] = 0
$$
and there are joint eigenfunctions of these commuting operators.
Such functions are called {\it Floquet} (or {\it Bloch}) functions.
The rigorous definition is as follows.

\begin{definition}
A function $f: \R^n \rightarrow \C$ is called a Floquet function
of a $\Lambda$-periodic operator $L$ if
$$
L f = E f
\ \ \ \ \mbox{and}\ \ \ \
f(x + \gamma) =
e^{2\pi i \langle k,\gamma \rangle} f(x)
$$
for $\gamma \in \Lambda$. The quantities $k_1,\dots,k_n$ are called
the quasimomenta of $f$.
\end{definition}

Any Floquet function defines the {\it multiplier} homomorphism
$\mu:\Lambda \to \C$:
$$
f(x + \gamma) = \mu(\gamma) \cdot f(x).
$$

Consider the case, when $L = {\cal D}$, the two-dimensional Dirac operator
$$
{\cal D} =
\left( \begin{array}{cc} 0 & \partial \\ -
\bar{\partial} & 0 \end{array} \right) +
\left(
\begin{array}{cc}
U & 0 \\
0 & U
\end{array}
\right)
$$
with a double-periodic continuous potential $U(z)$
where $z =x^1 + ix^2 \in \C$.

\begin{theorem}
There is an analytic set  $Q({\cal D}) \subset \C^3$ of positive codimension
such that there exists a Floquet function of ${\cal D}$ with the
quasimomenta $k= (k_1,k_2)$ and the eigenvalue $E$ if and only if
$(k_1,k_2,E) \in Q({\cal D})$.

Its intersection with the plane $\lambda = 0$ is an analytic set
$Q_0({\cal D)}$ of complex dimension one.
\end{theorem}

{\sl Proof.} Take a constant $C$ such that the operator ${\cal D}+C$ is
invertible on $L_2(T^2) = L_2(\C/\Lambda)$. Then consider a polynomial
operator pencil
$$
A_{k,E} =
1 +
\left(
\begin{array}{cc}
U - (C + E) & \pi(k_2 +  i k_1) \\
\pi(k_2 - i k_1) & U - (C + E)
\end{array}
\right)
\left(
\begin{array}{cc}
C & \partial \\
- \bar{\partial} & C
\end{array}
\right)^{-1}.
$$
Since for any function $g$ we have
$$
e^{2\pi i \langle k,x \rangle} [A_{k,E}({\cal D}+C)] g =
[{\cal D} - E] e^{2\pi i \langle k,x \rangle} g,
$$
there is a Floquet function $f(x)$ with the quasimomenta $(k_1,k_2)$ and the
eigenvalue $E$ if and only if there is a $\Lambda$-periodic function
$g(x)$ satisfying the equation
$$
A_{k,E} [{\cal D}+C] g = 0.
$$
Such a solution exists
if and only if there is
a solution $\varphi \in L_2(T^2)$ to the equation
\begin{equation}
A_{k,E}\varphi = 0.
\label{keldysh}
\end{equation}
If such a solution $\varphi$ exists, then
$f = e^{2\pi i \langle k,x \rangle}[{\cal D}+C]^{-1}\varphi$ is the desired
Floquet function.
The operator pencil
$$
1 - A_{k,E}
$$
is polynomial in $k_1,k_2$, and $E$.
Since $U$ is bounded, the multiplication operator
$$
\times U: L_2(T^2) \to L_2(T^2)
$$
is bounded and the pencil $(1 - A_{k,E})$
consists in compact operators
on $L_2(T^2)$.

Now we apply the Keldysh theorem (or the polynomial
Fredholm alternative) to it. This theorem reads that there is a
regularized determinant $\widetilde{\det}\, A_{k,E}$
of this pencil analytic in $k_1,k_2$, and $E$
such that the equation (\ref{keldysh}) is solvable if and only if
$\widetilde{\det}\, A_{k,E} = 0$ \cite{Kel,Kel2}. Now it remains to put
$$
Q({\cal D}) = \{\widetilde{\det}\, A_{k,E} = 0\}.
$$
In the same manner for the pencil $(1-A_{k,0})$ we obtain a complex
curve
$$
Q_0({\cal D}) = \{\widetilde{\det}\, A_{k,0} = 0\} \subset \C^2.
$$

As shown by perturbation methods the codimensions of these
sets are positive (see \cite{T3}). Also nontriviality of such determinants
follows from their construction by Keldysh.
This proves the theorem.

We applied this method to the operators
$\Delta+u$ and $\partial_t - \Delta$ in 1985. Later it
became known to us that Kuchment also proposed the same approach in \cite{Ku}
and therefore our paper was not published and referred in \cite{Kr,Kr2} as
an unpublished paper. The theory of such determinants is developed in
\cite{Ku2} and one can show that in fact they are entire
functions of $k$ and $E$.

We shall discuss another and more effective but technically difficult
approach of Krichever in 3.3.

Recall that the dual lattice $\Lambda^{\ast} \subset \C^2$
consists of vectors
$\gamma^{\ast}$ such that $\langle \gamma, \gamma^{\ast} \rangle =0$ for any
$\gamma \in \Lambda$.
The following proposition is evident.

\begin{proposition}
The sets $Q({\cal D})$ and $Q_0({\cal D})$ are invariant with respect to
translations by vectors from $\Lambda^{\ast}$:
\begin{equation}
k_1 \to k_1 + \Re \gamma^{\ast}, \ \ \
k_2 \to k_2 + \Im \gamma^{\ast}.
\label{dual}
\end{equation}
\end{proposition}

Indeed, this action preserves the multipliers.

In the sequel we shall confine to the set $Q_0({\cal D})$.

\begin{definition}
$Q_0({\cal D})$ is called the (zero)
Floquet spectral data of ${\cal D}$.
The genus of the normalization of $Q_0({\cal D})/\Lambda^{\ast}$ is called
the spectral genus of ${\cal D}$.
\end{definition}

Another two properties of the Floquet spectrum are easily derived in the
manner usual for soliton theory.

\begin{proposition}
If $U$ is real-valued, then
$Q_0({\cal D})$ is invariant under an antiholomorphic involution
$k \to -\bar{k}$.
\end{proposition}

{\sl Proof.} If $(\psi_1,\psi_2)^{\top}$ is a Floquet function
with the quasimomenta $(k_1,k_2)$, then
$(\bar{\psi}_2,-\bar{\psi}_1)^{\top}$ is a Floquet function with the
quasimomenta $(-\bar{k}_1,-\bar{k}_2)$.
This proves the proposition.

\begin{proposition}
$Q_0({\cal D})$ is invariant under a holomorphic involution
$k \to -k$.
\end{proposition}

{\sl Proof.}
Consider the pencil $L_k = A_{k,0}({\cal D}+C)$. We have
\begin{equation}
\left(
\begin{array}{cc}
0 & 1 \\ -1 & 0
\end{array}
\right) \cdot
\overline{L_k^{\ast}}
\cdot
\left(
\begin{array}{cc}
0 & -1 \\ 1 & 0
\end{array}
\right) =
L_{-k}.
\label{involution}
\end{equation}
The index of a Fredholm operator $A$ is
${\rm ind}\, A = \dim \ker A - \dim \ker A^{\ast}$.
Since ${\cal D}$ is selfadjoint, its index
vanishes. The operators $L_k$ have the same principal terms as ${\cal D}$
and, by the index theorem, their indices also vanish.
Hence if $k \in Q_0({\cal D})$, then $\dim L_k =
\dim L_k^{\ast} > 0$ and the identity
(\ref{involution}) implies that $\dim L_{-k}^{\ast} > 0$. Therefore
$(-k) \in Q_0({\cal D})$ and this proves the proposition.

Given a basis $(\gamma_1,\gamma_2)$ for $\Lambda$, we have a mapping
$$
{\cal M}: Q_0({\cal D})/\Lambda^{\ast} \to \C^2 \ \ \ :
\ \ \
{\cal M}(k) = (e^{2\pi i \langle k,\gamma_1 \rangle},
e^{2\pi i \langle k,\gamma_2 \rangle}),
$$
which maps quasimomenta into multipliers.

The submanifold ${\cal M}(Q_0({\cal D})/\Lambda^{\ast}) \subset \C^2$
is generically singular and its normalization is
the normalization of $Q_0({\cal D})/\Lambda^{\ast}$.

\begin{definition}
A complex curve $\Gamma$, which is
the normalization of $Q_0({\cal D})/\Lambda^{\ast}$,
is called the spectral curve of ${\cal D}$.
\end{definition}

A normalization sometimes consists in unstucking double points:
a pair of points of $\Gamma$ corresponding to a double
point of $Q_0({\cal D})/\Lambda^{\ast}$
is called a {\it resonance pair}.

The definition of ${\cal M}$ depends on a choice of a basis for
$\Lambda$. Given another basis $(\widetilde{\gamma}_1,\widetilde{\gamma}_2)$
for $\Lambda$ such that
$$
\left(
\begin{array}{c}
\widetilde{\gamma}_1 \\
\widetilde{\gamma}_2
\end{array}
\right) =
\left(
\begin{array}{cc}
a & b \\ c & d
\end{array}
\right)
\left(
\begin{array}{c}
\gamma_1 \\
\gamma_2
\end{array}
\right),  \ \ \ \
\left(
\begin{array}{cc}
a & b \\ c & d
\end{array}
\right) \in SL(2,\Z),
$$
the multipliers $(\mu_1,\mu_2) = (\mu(\gamma_1),\mu(\gamma_2))$ and
$(\widetilde{\mu}_1,\widetilde{\mu}_2) =
(\mu(\widetilde{\gamma}_1),\mu(\widetilde{\gamma}_2))$
are related as follows
\begin{equation}
\widetilde{\mu}_1 = \mu_1^a \mu_2^b, \ \ \
\widetilde{\mu}_2 = \mu_1^c \mu_2^d,
\label{sl2}
\end{equation}
and the sets of multipliers
$\{(\mu_1,\mu_2)\}$ and $\{(\widetilde{\mu}_1,\widetilde{\mu}_2)\}$
with respect to different bases are biholomorphically equivalent.

We call the image of ${\cal M}$ the {\it (Floquet) spectrum} of ${\cal D}$
(on the zero energy level, $E=0$). Given a basis of $\Lambda$,
this image is uniquely defined.
In general, the spectral data are defined modulo the
$SL(2,\Z)$-action (\ref{sl2})
and we say that the spectral data of two operators coincide if the
$SL(2,\Z)$-orbits of their spectral data coincide.

\subsection{Examples of spectra}

1) {\sl $U=0$.}

Let $\Lambda = \Z + i \Z$. The Floquet functions are
parameterized by two planes:
$$
\psi^+ = (e^{\lambda_+ z},0), \ \ \psi^- = (0,e^{\lambda_- \bar{z}}),
$$
$\Gamma$ is a union of these planes compactified by
two points at infinities, and $\psi$ has exponential
singularities at these points.

The quasimomenta of $\psi^+$ are
\begin{equation}
k_1 = \frac{\lambda_+}{2\pi i} + n_1, \ \
k_2 = \frac{\lambda_+}{2\pi} + n_2, \ \ n_j \in \Z,
\label{quasi1}
\end{equation}
and the quasimomenta of $\psi^-$ are
\begin{equation}
k_1 = \frac{\lambda_-}{2\pi i} + m_1, \ \
k_2 = -\frac{\lambda_-}{2\pi} + m_2, \ \ m_j \in \Z.
\label{quasi2}
\end{equation}
Hence
$$
Q_0 = \left(\cup_{n_1,n_2 \in \Z} A_{n_1,n_2}\right) \cup
\left(\cup_{m_1,m_2 \in \Z} B_{m_1,m_2}\right),
$$
where $A_{n_1,n_2}$ and $B_{m_1,m_2}$ are planes described by
(\ref{quasi1}) and (\ref{quasi2}).

The functions $\psi^+$ and $\psi^-$ have the same multipliers at the points
$$
\lambda_+^{m,n} = \pi(n+im), \ \ \lambda_-^{m,n} = \pi(n-im), \ \ m,n \in \Z.
$$
These are resonance pairs for this potential  with $\Lambda =
\Z + i \Z$.

Considering the zero potential as $\Lambda$-periodic with respect to
a general lattice $\gamma_1\Z + \gamma_2\Z$, the Floquet functions
are the same but the resonance pairs are
\begin{equation}
\lambda_+^{m,n} =
\frac{2\pi i}{\bar{\gamma}_1\gamma_2 - \gamma_1\bar{\gamma}_2}
(\bar{\gamma}_1 n - \bar{\gamma}_2 m), \ \
\lambda_-^{m,n} = \frac{2\pi i}{\bar{\gamma}_1\gamma_2 -
\gamma_1\bar{\gamma}_2}
(\gamma_1 n - \gamma_2 m).
\label{resonance}
\end{equation}

2) {\sl $U = C = \mbox{\rm const} \neq 0$.}

Assume for simplicity that $\Lambda =  \Z + i \Z$.

The Floquet functions are
$$
\psi(z,\bar{z},\lambda) = \left(
\exp{\left(\lambda z - \frac{C^2}{\lambda} \bar{z} \right)},
-\frac{C}{\lambda}\exp{\left(\lambda z - \frac{C^2}{\lambda} \bar{z} \right)}
\right)
$$
where $\lambda \in \Gamma = \C^{\ast} = \C \setminus \{0\}$.
Compactify $\Gamma$ by the points
$0$ and $\infty$ and define the Floquet function on $\Gamma$
globally as
$$
\psi(z,\bar{z},\lambda) = \frac{\lambda}{\lambda - C}
\left(\exp{\left(\lambda z - \frac{C^2}{\lambda} \bar{z} \right)},
- \frac{C}{\lambda}\exp{\left(\lambda z  - \frac{C^2}{\lambda}\bar{z}\right)}
\right).
$$
It has the following asymptotics
\begin{equation}
\psi \approx \left(\begin{array}{c} \exp{(k_+ z)} \\ 0 \end{array}\right)
\ \mbox{as $\lambda \to \infty$},\
\psi \approx \left(\begin{array}{c} 0 \\ \exp{(k_- \bar{z})} \end{array}\right)
\ \mbox{as $\lambda \to 0$}
\label{asymptotics}
\end{equation}
with $k_+ = \lambda$ and  $k_- = -C^2/\lambda$. After the
normalization $\psi$ gets a pole at $\lambda=C$.

The resonance pairs $(\lambda,\lambda')$ are
$$
\lambda = \frac{q\bar{q} \pm \sqrt{(q\bar{q})^2 - 4C^2 q\bar{q}}}{2\bar{q}},
\ \ \lambda' = \lambda - q, \ \ q=\pi(n+im), \ m,n \in \Z
$$
and they are parameterized by $q \in \pi \Z^2 \setminus \{0\}$.

3) {\sl $U = U(x)$ is a function of one variable.}

Let $U(x+T) = U(x)$ where $T$ is the minimal period.
Then the equation (\ref{diracequation}) for Floquet functions
$\psi(x,y) = \varphi(x)e^{\lambda y}$
is the {\it Zakharov--Shabat} system
\begin{equation}
L \varphi =
\left(\begin{array}{cc}
U & \frac{1}{2}\partial_x \\
- \frac{1}{2}\partial_x & U
\end{array}\right)
\left(\begin{array}{c}
\varphi_1 \\ \varphi_2
\end{array}\right) =
\left(\begin{array}{cc}
0 & \frac{i}{2}\lambda \\
\frac{i}{2}\lambda & 0
\end{array}\right)
\left(\begin{array}{c}
\varphi_1 \\ \varphi_2
\end{array}\right)
\label{zakharov}
\end{equation}
and in terms of $\eta_1 = \varphi_1 + i \varphi_2$ and
$\eta_2 = \varphi_1 - i\varphi_2$ it is
$$
(\partial_x + 2iU)\eta_1 = -i\lambda \eta_2, \ \
(\partial_x -2iU)\eta_2 = -i\lambda \eta_1.
$$
We see that $f = \eta_2$ satisfies the equation
$$
\left[-\partial_x^2 + (2 i U_x- 4U^2)\right] f = \nu f
$$
where $\nu = \lambda^2$. The transformation
$L \to \left[-\partial_x^2 + (2 i U_x- 4U^2)\right]$
is called the {\it Miura transformation}.
The same name is used for the transformation
\begin{equation}
\left[-\partial_x^2 + (2 i U_x- 4U^2)\right] \leftrightarrow
\left[-\partial_x^2 + (-2 i U_x- 4U^2)\right].
\label{miura}
\end{equation}

For any $\lambda \in \C$ take a two-dimensional space ${\cal V}_{\lambda}$
of solutions to (\ref{zakharov}) and consider the monodromy operator
$$
\widehat{T}:{\cal V}_{\lambda} \to {\cal V}_{\lambda}
\ \ : \ \
\widehat{T}(\varphi)(x) = \varphi(x+T).
$$
For any pair $(\varphi(x,\lambda),\vartheta(x,\lambda))$ of solutions to
(\ref{zakharov}) their Wronskian $W(\vartheta,\varphi) =
\vartheta_1(x) \varphi_2(x) - \vartheta_2(x) \varphi_1(x)$ is constant.
Take the basis $(c(x,\lambda),s(x,\lambda))$
for ${\cal V}_{\lambda}$ normalized as
$$
c(0,\lambda) = s_x(x,\lambda) = 1, \ \ c_x(0,\lambda) = s(x,\lambda) = 0.
$$
As shown in this basis the  entries of
the matrix $\widehat{T}$ are entire functions of $\lambda$.
Since  $W(c,s)$ is constant,
$\det \widehat{T}(\lambda) = 1$ and  the characteristic
equation for $\widehat{T}(\lambda)$ takes the form
\begin{equation}
k^2 - \Tr\widehat{T}(\lambda)\, k + 1 = 0.
\label{characteristic}
\end{equation}
If $\widehat{T}(\lambda)$ is diagonalized, then eigenfunctions
of $\widehat{T}$ are the Floquet functions.
The operator $\widehat{T}(\lambda)$ is not diagonalized if and only if
$\lambda$ is a simple root of the equation
\begin{equation}
\Tr^2\widehat{T}(\lambda) - 4 = 0,
\label{spur}
\end{equation}
which can have only simple and double roots.
In this case the Jordan form for $\widehat{T}(\lambda)$ is a non-diagonal
upper triangular matrix, and there is only one (up to multiple)
eigenfunction of $\widehat{T}(\lambda)$ with this value of $\lambda$.

We conclude that the Floquet function is globally defined on the
two-sheeted covering of $\C$ branched at points $\lambda_1,\dots$ which
are simple roots of (\ref{spur}) and this complex curve is exactly $\Gamma$.

Resonance pairs of the spectrum are pairs of points which project
into double roots of (\ref{spur}).

If there are finitely many
simple roots of (\ref{spur}) $L$ is called {\it finite gap}.
For finite gap operators $\Gamma$ is compactified by two infinities and
the Floquet functions are pasted in a meromorphic function on $\Gamma$
with the asymptotics (\ref{asymptotics}) at the infinities.

These analytic properties of the Floquet function for a one-dimensional
Schr\"o\-din\-ger operator are explained in \cite{DMN} and for the Dirac
operator such results are obtained by using the Miura transformation
or derived by the same reasonings straightforwardly.

\subsection{Spectra via perturbation theory}

Generically Floquet functions and the Floquet spectrum
can not be found by solving ordinary
differential equations as in examples in 3.2.
For proving existence of the Floquet spectrum we use in 3.1 the
Keldysh theorem. Another approach for finding this spectrum and
describing it in rather efficient manner was proposed by Krichever
who realized it for a two-dimensional Schr\"o\-din\-ger operator
and for the operator $\partial_y - \partial^2_x + U(x,y)$
\cite{Kr,Kr2}.  It is based on perturbation theory.

The examples discussed above demonstrate how the spectrum deforms under a
deformation of $U$. The main picture is as follows:
deforming potential we deform double points on
$Q_0({\cal D})/\Lambda^{\ast}$ into handles removing
singularities. The norm of the deformation measures the ``size'' of such
handles.

For the two-dimensional Dirac operator (\ref{diracoperator}) this
is not done until recently
but, since it is clear that the method of \cite{Kr}
works for this operator after a slight modification,
we explain what is the expected picture:

\begin{Pretheorem}
For a smooth potential $U$
the spectral curve $\Gamma$
consists of two parts: $M_0$ and $M_{\infty}$
where

1) $M_0$ is a complex curve of finite genus whose boundary consists in
a pair of circles;

2) $M_{\infty}$ is diffeomorphic to a union
of the domains $|\lambda_{\pm}| \geq R$ of the $\lambda_{\pm}$-planes
for some $R$
with  some resonance pairs $(\lambda_+^{m,n},\lambda_-^{m,n})$
(\ref{resonance}) ``joined by handles'' with decreasing sizes as
$m^2+n^2 \to \infty$;

3) $M_0$ and $M_{\infty}$ are pasted along their boundaries.

This ``joining by a handle'' means that some small disks
$|\lambda_+ - \lambda_+^{m,n}| < \varepsilon_+^{m,n}$ and
$|\lambda_- - \lambda_-^{m,n}| < \varepsilon_-^{m,n}$
are excluded and replaced by a cylinder pasted
to their boundaries and $\varepsilon_{\pm}^{m,n} \to 0$
as $m^2 + n^2 \to \infty$.

To any point of the subsets $M_+,M_- \subset M_{\infty}$, where
$$
M_+ = \C \setminus \{\ \{|\lambda_+| > R\} \cup \left(\cup_{m,n}
\{|\lambda_+ - \lambda_+^{m,n}| \leq  \varepsilon_+^{m,n}\}\right)\ \},
$$
$$
M_- = \C \setminus \{\ \{|\lambda_-| > R\} \cup \left(\cup_{m,n}
\{|\lambda_- - \lambda_-^{m,n}| \leq  \varepsilon_-^{m,n}\}\right)\ \},
$$
corresponds a unique (up to multiple) Floquet function.
These functions and
their multipliers $\mu(\gamma_j,\lambda_{\pm})$ asymptotically
behave as in the case $U = 0$, and, in particular,
$$
\mu(\gamma_j,\lambda_+) = e^{\lambda_+ \gamma_j}
\left(1 + O\left(\frac{1}{\lambda_+}\right)\right), \ \
\mu(\gamma_j,\lambda_-) = e^{\lambda_- \bar{\gamma}_j}
\left(1 + O\left(\frac{1}{\lambda_-}\right)\right)\ \
$$
as $\lambda_{\pm} \to \infty$.

If $U$ does not vanish identically, then $\Gamma$ is irreducible.
\end{Pretheorem}

A potential $U$ is finite gap if there is such a
representation with no handles joining resonance points in
$M_{\infty}$. In this case $\Gamma$ is compactified by a pair of
infinities $\infty_{\pm}$ to a Riemann surface of finite genus.

In fact this is a general description of the Floquet spectra of operators
with periodic coefficients.
This physical picture from \cite{Kr} was chosen
by Feldman, Kn\"orrer, and Trubowitz as a most convenient and reasonable
definition of general (non-hyperelliptic) Riemann surfaces of infinite
genus and they developed a nice theory of such
surfaces for which analogs of many classical theorems on algebraic curves
take place \cite{FKT} (see also their preprints published by ETH).

\section{The spectrum of the Weierstrass representation}

\subsection{The spectral curve of an immersed torus}

Let $F: \C \to \R^3$ be a conformal immersion of a plane whose
Gauss map descends through a torus
$\C/\Lambda$, i. e., double periodic.
We shall consider tori as a particular case of such planes,
when the immersion is also double periodic.

For immersed plane with a periodic Gauss map we have the
Weierstrass representation constructed in 2.2.
The double periodic
potential $U(z)$ of this representation is the potential of the surface
(with fixed conformal parameter).

It is said that two surfaces (with fixed conformal parameters)
$F_1:\C \to \R^3$ and $F_2:\C \to \R^3$
are {\it isopotential} if their potentials coincide.

\begin{definition}
The spectral curve $\Gamma$ of the operator ${\cal D}$ with
the potential $U$ is called the {\it spectral curve of a surface}
and the spectral genus of ${\cal D}$ is called
the spectral genus of a surface. Given a basis
$\gamma_1, \gamma_2$ for $\Lambda$, the image
of the multiplier map
$$
{\cal M}: Q_0({\cal D})/\Lambda^{\ast} \to \C^2 \ \ \ : \ \ \
{\cal M}(k) = (e^{2\pi i \langle k,\gamma_1 \rangle},
e^{2\pi i \langle k,\gamma_2 \rangle})
$$
is called the spectrum of a surface.
\end{definition}

Let us look how these spectral notions
depend on a choice of a conformal parameter.

\begin{proposition}
The spectral curve and the spectral genus of a surface do not depend
on a choice of a conformal parameter.

The spectrum of a surface
depends on a choice of a basis for $H_1(T^2) = H_1(\C/\Lambda)= \Lambda$
and for different bases $(\widetilde{\gamma}_1,\widetilde{\gamma}_2)$ and
$(\gamma_1,\gamma_2)$ related by a $SL(2,\Z)$-transformation
$$
\left(
\begin{array}{c}
\widetilde{\gamma}_1 \\
\widetilde{\gamma}_2
\end{array}
\right) =
\left(
\begin{array}{cc}
a & b \\ c & d
\end{array}
\right)
\left(
\begin{array}{c}
\gamma_1 \\
\gamma_2
\end{array}
\right),
$$
the spectra $\{(\widetilde{\mu}_1,\widetilde{\mu}_2)\}$ and
$\{(\mu_1,\mu_2)\}$
are related as
$$
\widetilde{\mu}_1 = \mu_1^a \mu_2^b, \ \ \
\widetilde{\mu}_2 = \mu_1^c \mu_2^d.
$$
\end{proposition}

{\sl Proof.}
It remains to check that all these data are preserved by
a passage from $z$ to $w = t^2z$ with $t \in \C \setminus \{0\}$.
For these conformal parameters the surface is defined by
functions $\psi(z) = (\psi_1(z),\psi_2(z))^\top$ and
$\widetilde{\psi}(w) = (\widetilde{\psi}_1(w),\widetilde{\psi}_2(w))^\top$,
which are solutions to the equations
$$
{\cal D}\psi =
\left(
\begin{array}{cc}
U & \partial_z \\ -\partial_{\bar{z}} & U
\end{array}
\right)\psi = 0, \ \
\widetilde{{\cal D}}\widetilde{\psi} =
\left(
\begin{array}{cc}
\widetilde{U} & \partial_w \\ -\partial_{\bar{w}} & \widetilde{U}
\end{array}
\right)\widetilde{\psi} = 0.
$$
We have
$d w = t^2 dz, \partial_z = \frac{1}{t^2} \partial_w$,
and, since
$e^{2\widetilde{\alpha}(w)}d w d\bar{w} =
e^{2\alpha(z)} d z d \bar{z}$ and $\widetilde{H}(w) = H(z)$,
the formulas (\ref{formulas}) imply that
$\widetilde{U}(w) = t\bar{t}\, U(z)$ for $w = t^2z$.
Therefore
$$
\left(
\begin{array}{cc}
U & \partial_z \\ -\partial_{\bar{z}} & U
\end{array}
\right)
\left(
\begin{array}{c}
\varphi_1 \\ \varphi_2
\end{array}
\right) = 0
\ \ \
\mbox{if and only if}
\ \ \
\left(
\begin{array}{cc}
\widetilde{U} & \partial_w \\ -\partial_{\bar{w}} & \widetilde{U}
\end{array}
\right)
\left(
\begin{array}{c}
t\,\varphi_1 \\ \bar{t}\,\varphi_2
\end{array}
\right) = 0
$$
and, since $\varphi=(\varphi_1,\varphi_2)^\top$ and
$\widetilde{\varphi} = (t\varphi_1,\bar{t}\varphi_2)^\top$
have the same multipliers.

The transformation of the spectra was already discussed in 3.1 and we
just repeat here these formulas because now they appear in another
situation.

The proposition is proved.

We say that two planes, which may convert to tori by immersions,
are {\it isospectral} if

1) there are conformal parameters one them such
that both of them are represented by mappings
$F_1:\C \to \R^3$ and $F_2:\C \to \R^3$;

2) the corresponding potentials
of Dirac operators are $\Lambda$-periodic with the same lattice $\Lambda$,
and the spectra of these operators
with respect to a fixed basis for $\Lambda$ coincide.

\subsection{On the Willmore functional}

Given a Weierstrass representation of a torus $\Sigma$, we have
\begin{equation}
{\cal W}(\Sigma) =
4 \int_{\Pi} U^2 d x \wedge d y = 2i\int_{\Pi} U^2 dz \wedge d\bar{z},
\label{willmore}
\end{equation}
(see \cite{T1}) and that shows that
the Willmore functional measures the deviation of ${\cal D}$ from
the Dirac operator with the zero potential and geometrically that means that
it measures the deviation of a connection in a spinor bundle defined by
${\cal D}$ from the trivial connection.

A relation of this spectrum to the Willmore functional
$$
{\cal W}(\Sigma) = \int_{\Sigma} H^2 \, d\,\mu,
$$
was first established in
\cite{T2} where it was discussed for surfaces of revolution.
In this case there is a direct construction of the Floquet spectrum
which a hyperelliptic complex curve (see Example 3 in 3.2).
When the curve is of finite genus there is a compactification of it
by two infinities and the Floquet function $\psi(z,P)$
is defined on the compactification and is meromorphic outside these infinities.
We shall discuss a general case using
a construction of the spectrum by perturbation.

Let $U$ be a $\Lambda$-periodic potential with
$\Lambda = \gamma_1\Z + \gamma_2\Z$ and let
$\Gamma$ be of finite genus and be
compactified by two infinities $\infty_{\pm}$. In fact, as shown in
Pretheorem,  these infinities are
inherited during the perturbation of $U$
from the compactification of the spectrum of the zero potential
(see Example 1 in 3.2).

Take the Floquet function $\psi(z,P)$
meromorphic outside the infinities
and with the asymptotics
\begin{equation}
\psi \approx \left(\begin{array}{c} \exp{(\lambda_+ z)} \\ 0
\end{array}\right)
\ \mbox{as $P \to \infty_{\pm}$},\
\psi \approx \left(\begin{array}{c} 0 \\ \exp{(\lambda_- \bar{z})}
\end{array}\right)
\ \mbox{as $P \to \infty_-$}
\label{asymptotics2}
\end{equation}
where $\lambda^{-1}_{\pm}$ are local parameters near $\infty_{\pm}$.

The theory of finite gap integration gives a recipe for
reconstructing $U$ from such asymptotic expansions.
Let
$$
\psi(z,\lambda_+) = \exp{(\lambda_+ z)} \cdot
\left(
\left(\begin{array}{c} 1 \\ 0 \end{array}\right) +
\left(\begin{array}{c} \zeta_1 \\ \zeta_2 \end{array}\right)
\frac{1}{\lambda_+} +
O\left(\frac{1}{\lambda_+^2}\right)\right)
\
\mbox{as $\lambda_+ \to \infty$.}
$$
Substitute this expansion into the Dirac equation and
expand ${\cal D}\psi = 0$ into the powers of $\lambda_+$.
Every coefficient in this expansion equals zero. Take the first two of them:
\begin{equation}
U + \zeta_2 = 0, \ \ \  U\zeta_2 - \bar{\partial}\zeta_1 = 0.
\label{reconstruction}
\end{equation}
This gives a reconstruction formula for $U$ and also the identity
$$
- U^2 = \bar{\partial} \zeta_1.
$$

Now let us show how the Willmore functional appears in this picture.
By the perturbation theory, we expect that the spectrum
asymptotically behaves as the spectrum of the zero potential and this
leads to the following conclusion: there is a function $W(\lambda_+)$
defined near $\infty_+$ such that

1) $W(\lambda_+) = C_1 \lambda_+^{-1} +
O\left(\lambda_+^{-2}\right)$;

2) $\psi(z,\lambda_+) = e^{\lambda_+z + W(\lambda_+)\bar{z}}
\cdot \varphi(z,\lambda_+)$ where $\varphi(z,\lambda_+)$ is
$\Lambda$-periodic.

The function $W$ measures the deviation of the Floquet spectrum from the
spectrum of the zero potential.
Indeed, the multipliers of $\psi$ are
$$
(\mu(\gamma_1),\mu(\gamma_2)) = (e^{\lambda_+ \gamma_1 +
W(\lambda_+)\bar{\gamma}_1},
e^{\lambda_+ \gamma_2 + W(\lambda_+)\bar{\gamma}_2})
$$
and the multipliers of the Floquet function of the zero potential
(which is considered as $\Lambda$-periodic)  are
$$
(\mu_0(\gamma_1),\mu_0(\gamma_2)) =
(e^{\lambda_+ \gamma_1},e^{\lambda_+ \gamma_2}).
$$
It is easy to notice that
$$
\zeta_1 = C_1 \bar{z} \ + \ \mbox{( a $\Lambda$-periodic function)}
$$
and hence
$$
\int_{\Pi} U^2 \, d z \wedge d\bar{z} = - C_1 \cdot \mbox{\rm Vol}\,\Pi
$$
where $\Pi$ is a the
parallelogram spanned by $\gamma_1$ and $\gamma_2$.
Now using (\ref{willmore}) we conclude that
\begin{equation}
{\cal W}(\Sigma) = 4 C_1 \cdot (\gamma_1 \bar{\gamma}_2 -
\bar{\gamma}_1 \gamma_2) = - 4 C_1 \cdot \mbox{\rm Vol}\,\Pi
\label{residue}
\end{equation}
where $\Pi$ is the area of the fundamental domain of $\Lambda$.
This derivation of the formula (\ref{residue}) was exposed
in \cite{GS}. The analogous derivation for the area of minimal tori is given
in \cite{Hitchin}.

Now consider the whole series
$$
W(\lambda) = C_1\lambda^{-1} + C_2\lambda^{-2} + \dots.
$$
Since an involution $k \to -k$ inverts $\lambda$ and
preserves the spectrum, we have $C_{2k} = 0$ for $k=1,2,\dots$.
The quantities $C_{2k-1}$ for $k \geq 2$
depend on choices of a conformal parameter $z$
and a parameter $\lambda$ on the spectral curve.
Given a parameter $\lambda$,
the holomorphic differentials
$$
{\cal W}_k = - 4 (C_{2k-1}\, \mbox{\rm Vol}\,\Pi) \, dz^{2k-2}
$$
are geometric invariants.
The first of them is the Willmore functional ${\cal W} = {\cal W}_1$.

If the conformal parameter is fixed, then $C_{2k-1}$
are first integrals of the modified Novikov--Veselov equation and
the question of what are their geometrical meanings was posed in \cite{T1}.

On a surface of revolution there is a distinguished conformal parameter
$z = x+ i y$ where $x$ is a parameter on the rotating curve and
$y$ is the angle of rotation, $0 \leq y \leq 2\pi$, and there is
a distinguished parameter $\lambda$ which is the eigenvalue of
the Miura transformation of the Dirac operator
(see 3.2). For these parameters
$C_{2k-1}$ are the {\it Kruskal--Miura integrals} of the
modified Korteweg--de Vries equations. They are
geometric invariants of surfaces of revolution \cite{T2}.

Critical points of the Willmore functional are {\it Willmore surfaces}.
For higher functionals critical points were not studied but
using trace formulas  we showed in \cite{T2} that for spheres of revolution
higher invariants are not bounded both from above and below.

\subsection{Surfaces in terms of theta functions}

Assume that the spectrum $\Gamma$ of a torus $\Sigma$ has finite genus
which equals $g$ and take two different points $\infty_{\pm}$ on $\Gamma$
with local parameters $\lambda^{-1}_{\pm}$ near
them such that $\lambda^{-1}_{\pm}
(\infty_{\pm})=0$.
Then the theory of Baker--Akhieser functions \cite{Kr00}
(see also \cite{DKN}) reads that for a generic effective
divisor $D$, i. e. , a formal sum of points on $\Gamma$, of degree
$g+1$ ($D = P_1 + \dots + P_{g+1}$) there is a unique function
$\psi(z,\bar{z},P)$ such that

1) $\psi$ is meromorphic in $P \in \Gamma \setminus \{\infty_{\pm}\}$
and has the asymptotics (\ref{asymptotics2}) as $P \to \infty_{\pm}$;

2) $\psi$ has poles only in $D$.

Then this function is constructed in terms of theta functions of
a Riemann surface $\Gamma$. From this function one can
reconstruct a Dirac operator ${\cal D}$ by (\ref{reconstruction}).
Therefore to each point $P \in \Gamma \setminus\{\infty_{\pm},
P_1,\dots, P_{g+1}\}$ there corresponds a surface in $\R^3$ constructed
from $\psi(P)$ via (\ref{weierstrass}).

In this event, the Riemann surface $\Gamma$ parameterizes isopotential
surfaces and each of these surfaces is described in terms of
theta functions of $\Gamma$.
The detailed formulas are given in \cite{T3,T5}. For CMC tori in $\R^3$
such formulas were derived in \cite{Bo1}.

For tori of infinite spectral genera one can apply the theory of theta
functions on such surfaces developed by
Feldman, Kn\"orrer, and Trubowitz \cite{FKT}.

\subsection{On Willmore surfaces and the Willmore conjecture}

First the Willmore functional ${\cal W}$
appeared in the 20s in the papers by
Blaschke \cite{Bl} and Thomsen \cite{Th}. In this event it was called
the {\it conformal area} and its extrema were called {\it conformally minimal
surfaces}. Blaschke and Thomsen also established the main properties of this
functional which are:

1) the Willmore functional is invariant with respect to conformal
transformations of the ambient space: let $\Sigma \subset \R^3$ be a
compact immersed oriented surface, $z = x+iy$ be a conformal parameter on
it, and $G: \bar{\R}^3 \to \bar{\R}^3$ be a conformal transformation which
maps $\Sigma$ into $\R^3$, then ${\cal W}(\Sigma) = {\cal W}(G(\Sigma))$,
and this follows from the conformal invariance of the quantity
\footnote{Pinkall indicated that he and Pedit
recently proved that the quantity $Ae^{-\alpha}$
is already conformally invariant.}
$$
(k_1-k_2)^2 d\mu = 4(H^2 - K)d\mu = 16|A|^2 e^{-2\alpha} dx \wedge dy
$$
where $e^{2\alpha}dzd\bar{z}$ is the first fundamental form and $Adz^2$
is the Hopf differential of the surface;

2) if $\Sigma$ is a minimal surface in $S^3$ and $\pi:S^3 \to \bar{\R}^3$
is the stereographic projection which maps $\Sigma$ into $\R^3$, then
$\pi(\Sigma)$ is a conformally minimal (Willmore) surface.

Hence
in difference with minimal surfaces there are compact
immersed Willmore surfaces in $\R^3$.

All Willmore spheres were described by Bryant \cite{Bryant}.
The classification of Willmore tori is not complete until recently and we
only mention the papers \cite{Babich,BB}
where the finite gap integration was applied to this problem.
In \cite{Helein} the Dorfmeister--Pedit--Wu (DPW) method \cite{DPW},
was applied for constructing general Willmore surfaces.

The simplest example of a Willmore torus is the stereographic
projection of the Clifford torus $\{(x^1)^2+(x^2)^2=1/2, (x^3)^2+(x^4)^2 =
1/2\} \subset S^3 \subset \R^4$ into $\R^3$. In another way it may be obtained
as a circle torus of revolution such that the ratio of the distance from the
center of the circle to the axis of revolution and the radius of the circle
equals $\sqrt{2}$. This torus in $\R^3$ is also called the Clifford torus.

Willmore conjectured that

{\sl the Willmore functional
achieves its minimum for tori on the Clifford torus and its conformal
transformations and this minimum equals $2\pi^2$
(the Willmore conjecture)}

\noindent
and checked this conjecture for circle tori of
revolution \cite{Willmore}.

This conjecture implies the Hsiang--Lawson conjecture that the area of a
minimal torus in $S^3$ is greater or equal than $2\pi^2$, the area of the
Clifford torus in $S^3$, but not equivalent to it since there are
Willmore tori in $\R^3$
which are not images of minimal tori under the stereographic projection
\cite{Pinkall0}.

The Willmore conjecture is still open and there are some particular cases
for which it was proved:

1) Langer and Singer proved it for tori of revolution \cite{LS} and
Hertrich-Jeromin and Pinkall generalized their result for channel tori
which are the enveloping tori for one-parameter families of spheres
\cite{HJP};

2) Li and Yau proved it for tori
conformally equivalent to flat tori $\C/\{\Z + \tau \Z\}$ where
$\tau = a+ib, 0 \leq a \leq 1/2, b >0$ and $\sqrt{1-a^2} \leq b \leq 1$
\cite{LY} and later Montiel and Ros improved the latter inequality to
$(a-1/2)^2 + (b-1)^2 \leq 1/4$ \cite{MR}.

Simon proved that the minimum of the Willmore functional is achieved
on a real-analytic torus \cite{Simon}.

For higher genera surfaces a generalization of the Willmore conjecture was
proposed by Kusner \cite{Kusner}. He conjectured that for such surfaces the
Willmore functional attain its minima of the stereographic projections of
some minimal surfaces in $S^3$ constructed by Lawson.

In \cite{T3} (see also \cite{T5}) we conjectured that

{\sl for a fixed conformal classes of tori the minimum of ${\cal W}$ is
attained on tori with minimal spectral genus.}

If this conjecture is valid then we may reduce the Willmore conjecture to
estimating ${\cal W}$ for Willmore tori of small spectral genera and this
can be done by using soliton theory.

The global Weierstrass representation
gives a physical explanation for lower bounds for ${\cal W}$:
for small perturbations of the zero-potential
$U=0$ the surfaces constructed from solutions to
(\ref{diracequation}) via (\ref{weierstrass})
do not convert into tori and the lower bound for ${\cal W}$, the squared
$L_2$-norm of the perturbation, shows how large a perturbation has to be
to convert planes into tori.

\section{The spectra of integrable tori}

\subsection{Constant mean curvature tori}

\medskip

Let $\Sigma$ be a CMC torus in $\R^3$, i.e.,
$H = \mbox{const}$, and let it be
conformally equivalent to $\C/\Lambda$ with
$\Lambda = \gamma_1 \Z + \gamma_2 \Z$.

As shown in 2.3 for CMC surfaces
the Hopf differential $\omega = A\, d z^2$ is holomorphic
and, since
the space of quadratic holomorphic differentials on a torus
is one-dimensional, we have $A \,d z^2 = \mbox{const}\,d z^2$.
This differential does not vanish because otherwise all points are
umbilics which is impossible for tori in $\R^3$.
Hence assume that
$$
\omega = \frac{1}{2}\, d z^2, \ \ H=1
$$
and this is achieved by rescaling $z$ and by a homothety in $\R^3$.
Now the Codazzi equations in terms of $u=2\alpha$ read
\begin{equation}
u_{z\bar{z}} + \sinh u = 0
\label{sinh}
\end{equation}
which is the {\it sinh-Gordon equation}.
The Codazzi equations (\ref{w-codazzi})
give a commutation representation for (\ref{sinh}):
\begin{equation}
\left[
\frac{\partial}{\partial z} -
\left(
\begin{array}{cc}
\alpha_z & - \frac{\lambda^2}{2} e^{-\alpha} \\
- \frac{1}{2}e^{\alpha} & 0
\end{array}
\right)\right]\psi = 0,
\ \ \ \
\left[
\frac{\partial}{\partial \bar{z}} -
\left(
\begin{array}{cc}
0 & \frac{1}{2}e^{\alpha} \\
\frac{\lambda^{-2}}{2} e^{-\alpha} & \alpha_{\bar{z}}
\end{array}
\right)\right] \psi = 0
\label{sinh-comm-1}
\end{equation}
where $\lambda^2 = - 1$ in (\ref{w-codazzi}).

Notice that

1) for any $\lambda \neq 0$ the compatibility condition for
(\ref{sinh-comm-1}) is (\ref{sinh});

2) if $|\lambda|=1$ then (\ref{sinh-comm-1}) are the Codazzi equations
(\ref{w-codazzi}) for the surface defined by $\psi(\lambda,z,\bar{z})$
via (\ref{weierstrass}).

There is another representation of (\ref{sinh}) which gives rise to
the {\it spectral curve of a CMC torus}
\cite{PS,Bo1}.

Consider the equation (\ref{sinh}) as the compatibility condition for the
system
\begin{equation}
\left[
\frac{\partial}{\partial z} -
\frac{1}{2}
\left(
\begin{array}{cc}
-u_z & - \lambda \\
- \lambda & u_z
\end{array}
\right)\right]\varphi = 0, \ \ \
\left[
\frac{\partial}{\partial \bar{z}} -
\frac{1}{2 \lambda}
\left(
\begin{array}{cc}
0 & e^{-u} \\
e^u & 0
\end{array}
\right)\right] \varphi = 0
\label{sinh-comm-2}
\end{equation}
which contains the linear problem
$$
L \varphi =
\partial_z \varphi -
\frac{1}{2}
\left(
\begin{array}{cc}
-u_z & 0 \\
0 & u_z
\end{array}
\right) \varphi =
\frac{1}{2}
\left(
\begin{array}{cc}
0 & - \lambda \\
- \lambda & 0
\end{array}
\right) \varphi
$$
for a general $\Lambda$-periodic potential $u$.
Since $L$ is a first order $2\times 2$-matrix operator,
for every $\lambda \in \C$ the system
(\ref{sinh-comm-2}) has a two-dimensional space ${\cal V}_{\lambda}$
of solutions and these spaces are invariant under the translation operators
$$
\widehat{T}_j \varphi (z) = \varphi (z + \gamma_j), \ \ \ j=1,2.
$$
Since $\widehat{T}_1, \widehat{T}_2$, and $L$ commute,
they have common eigenvectors and these vectors
are glued into  a meromorphic function $\Psi (z,\bar{z},P)$ on a
two-sheeted covering
$$
\Gamma(L) \rightarrow \C: P \in \Gamma(L) \rightarrow \lambda \in \C,
$$
which ramifies at points where $\widehat{T}_j$ and $L$ are not
diagonalized simultaneously.

To each point $P \in \Gamma(L)$
corresponds a unique (up to a constant multiple)
Floquet function $\varphi$ with multipliers $\mu(\gamma_1,P)$
and $\mu(\gamma_2,P)$.

By the same reasonings as in the example 3 in 3.4 it is shown that
there are the Floquet functions defined on a $\Gamma(L)$
such that $\Gamma(L)$ is compactified by four infinities $\infty_{\pm}^1,
\infty_{\pm}^2$ such that
$\infty_{\pm}^1$ are mapped into $\lambda = \infty$ and $\infty_{\pm}^2$
are mapped into $\lambda = 0$ and there are asymptotics
$$
\psi(z,P) \approx
\exp{\left(\mp\frac{\lambda z}{2}\right)}
\left(
\begin{array}{c} 1 \\ \pm 1
\end{array}\right)
\ \ \mbox{as $P \to \infty_{\pm}^1$},
$$
$$
\psi(z,P) \approx
\exp{\left(\mp \frac{\bar{z}}{2\lambda}\right)}
\left(
\begin{array}{c} 1 \\ \pm 1
\end{array}\right)
\ \ \mbox{as $P \to \infty_{\pm}^2$},
$$
and therefore their multipliers tend to $\infty$ as $\lambda \to 0,\infty$
\cite{Bo1}.

If $\varphi = (\varphi_1,\varphi_2)^\top$ satisfies
(\ref{sinh-comm-2}) for $\lambda = \mu$, then
\begin{equation}
\sigma(\varphi) = (\varphi_1,-\varphi_2)
\label{sigma-cmc}
\end{equation}
satisfies (\ref{sinh-comm-2}) for $\lambda = -\mu$ and this generates
an involution $\sigma:\Gamma(L) \to \Gamma(L)$,
descending to an involution of $\C:\lambda \to -\lambda$.
We also have
$\sigma(\infty^1_{\pm}) = \infty^1_{\mp}$ and
$\sigma(\infty^2_{\pm}) = \infty^2_{\mp}$.
By (\ref{sigma-cmc}), the immersion
$$
{\cal M}:\Gamma(L) \to \C^2 \ \ : P \to (\mu(\gamma_1,P),\mu(\gamma_2,P))
$$
descends through the quotient of $\sigma$, i.e.,
${\cal M}:\Gamma(L) \to \Gamma(L)/\sigma \to \C^2$.

\begin{definition}
The complex curve $\Gamma(L)/\sigma$
is called the {\sl spectral curve}
of a CMC torus $\Sigma$.
It is said that ${\cal M}(\Gamma(L)/\sigma)$ is the spectrum of this torus.
\end{definition}

By straightforward computations, we obtain

\begin{proposition}
$\varphi=(\varphi_1,\varphi_2)^{\top}$ satisfies (\ref{sinh-comm-2}) if
and only if $\psi = (\lambda \varphi_2, e^{\alpha}\varphi_1)^{\top}$
satisfies (\ref{sinh-comm-1}).
\end{proposition}

This proposition implies that the Floquet functions
of $L$ and ${\cal D}$ have the same multipliers:

\begin{theorem}
Given a CMC torus $\Sigma$, the spectrum of the CMC torus
form a component, of the spectrum of the surface as defined in 4.1,
containing both
asymptotic ends where $\mu(\gamma_j) \approx e^{\lambda_+ \gamma_j}$ and
$\mu(\gamma_j) \approx e^{\lambda_- \bar{\gamma}_j}$ as $\lambda_{\pm}
\to \infty$.

Therefore, the spectral curve of the CMC torus $\Sigma$ is
an irreducible component of the spectral curve of this surface as defined
in 4.1.
\end{theorem}

Assuming that Pretheorem is valid, we have more:

{\sl the spectrum and the spectral curve of a CMC torus coincide
with the spectrum and the spectral curve of this surface as defined in 4.1.}

In fact, for this conclusion we need only one part of Pretheorem, which
states that the spectral curve is irreducible for $U \neq 0$.

Notice that the genus of the
spectral curve of a CMC torus is finite \cite{PS}.

\subsection{Isothermic tori}

A surface is called {\sl isothermic} if
near every point there is a conformal parameter
$z = x+iy$ such that the fundamental forms are
$$
{\bf I} = e^{2\alpha}(d x^2 + d y^2), \ \
{\bf II} = e^{2\alpha}(k_1 d x^2 + k_2 d y^2).
$$
To each isothermic surface $F:{\cal U} \to \R^3$ there corresponds
the {\it dual isothermic surface} $F^{\ast}:{\cal U} \to \R^3$ defined up to
translations by the formulas
$$
F^{\ast}_z = e^{-2\alpha}F_{\bar{z}}, \  \
F^{\ast}_{\bar{z}} = e^{-2\alpha}F_{z}.
$$
The fundamental forms of the dual surface are
$$
{\bf I} = e^{-2\alpha}(d x^2 + d y^2), \ \
{\bf II} = - k_1 d x^2 + k_2 d y^2,
$$
the Gauss maps of $F$ and $F^{\ast}$ are antipodal: $N = -N^{\ast}$,
and $F^{\ast\ast} = F$ (modulo translations).

It is obtained by straightforward computations that

\begin{proposition}
\label{Pisothermic}
If an isothermic surface $F$ is represented via (\ref{weierstrass}) by
a vector function $\psi = (\psi_1,\psi_2)^{\top}$, then the
dual surface $F^{\ast}$ is represented via (\ref{weierstrass})
by the function $\psi^{\ast} = (i e^{-\alpha}\psi_2,i e^{-\alpha}\psi_1)$.

The potentials of these surfaces are
$$
U = \frac{k_1+k_2}{4} e^{\alpha}; \ \
U^{\ast} = \frac{k_2- k_1}{4} e^{\alpha}
$$
and the Hopf differentials $A\, d z^2$ and $A^{\ast}\, d z^2$ are
$$
A = \frac{k_1-k_2}{4} e^{2\alpha}, \ \ \
A^{\ast} = -\frac{k_1+k_2}{4}.
$$
\end{proposition}

\begin{corollary}
\label{dualpotential}
Given an isothermic surface $\Sigma$ and its Hopf differential $A dz^2$
and the metric $e^{2\alpha}d z d\bar{z}$, there is an equality
$$
A e^{-\alpha} = \frac{k_1-k_2}{4} e^{\alpha} = -U^{\ast},
$$
where $U^{\ast}$ is the potential of the dual surface.
\end{corollary}

The simplest examples of isothermic surfaces are
surfaces of revolution and constant mean curvature surfaces.

Let $\Sigma$ be an isothermic plane, which may convert into an immersed
torus, whose Gauss map descends through
$\C/\Lambda$ with $\Lambda = \gamma_1\Z + \gamma_2\Z$.

By Proposition \ref{Pisothermic},
the Gauss--Weingarten equations (\ref{w-gw}) are written as
$$
\psi_z = {\bf U} \psi,  \ \  \psi_{\bar{z}} = {\bf V} \psi,
$$
with
$$
{\bf U} =
\left(\begin{array}{cc}
 \alpha_z & \frac{k_1 - k_2}{4} e^{\alpha} \\
- \frac{k_1 + k_2}{4} e^{\alpha} & 0
\end{array}
\right), \ \ \
{\bf V} =
\left(\begin{array}{cc}
0 & \frac{k_1 + k_2}{4} e^{\alpha} \\
\frac{k_2 - k_1}{4} e^{\alpha} & \alpha_{\bar{z}}
\end{array}
\right),
$$
and their compatibility conditions are
\begin{equation}
\alpha_{x x} + \alpha_{y y} +
k_1 k_2 e^{2\alpha} = 0, \ \ \
k_{2x} - (k_1 - k_2) \alpha_x = k_{1y} + (k_1 - k_2) \alpha_y = 0.
\label{isothermic}
\end{equation}

The equations (\ref{isothermic}) are also the compatibility
conditions for linear problems with a spectral parameter:
\begin{equation}
\widehat{\varphi}_z = \widehat{\bf U}(\lambda) \widehat{\varphi}, \ \ \
\widehat{\varphi}_{\bar{z}} = \widehat{\bf V}(\lambda) \widehat{\varphi},
\label{la-isothermic}
\end{equation}
where
$$
\widehat{\bf U}(\lambda) =
\left(
\begin{array}{cc}
{\bf U} & \lambda {\bf J}^- \\
\lambda {\bf J}^+ & {\bf U} + \alpha_z {\bf E}
\end{array}
\right), \ \
\widehat{\bf V}(\lambda) =
\left(
\begin{array}{cc}
{\bf V} & \lambda {\bf J}^+ \\
\lambda {\bf J}^- & {\bf V} + \alpha_{\bar{z}} {\bf E}
\end{array}
\right), \ \
$$
$$
{\bf J}^+ =
\left(
\begin{array}{cc}
0 & 1 \\ 0 & 0
\end{array} \right), \ \
{\bf J}^- =
\left(
\begin{array}{cc}
0 & 0 \\ 1 & 0
\end{array} \right), \ \
{\bf E} =
\left(
\begin{array}{cc}
1 & 0 \\ 0 & 1
\end{array} \right).
$$
First such representation with a spectral parameter had been found in
\cite{CGS} in terms of $5 \times 5$-matrices and, by using
the $4$-dimensional spinor representation of $SO(5,\C)$,
had been written in terms of $4 \times 4$-matrices in \cite{BP}.
Here we use another representation which is gauge equivalent to
the latter one.

Now the reasonings for describing the spectrum of an isothermic torus
are the same as for a one-dimensional Dirac operator
(see Example 3 in 3.2) and CMC tori and we will only sketch them.

For every $\lambda \in \C$ the system (\ref{la-isothermic}) has
a four-dimensional space ${\cal V}_{\lambda}$ of solutions.
On each such a space the translation operators act
$$
\widehat{T}_j \widehat{\varphi} (z) =
\widehat{\varphi} (z + \gamma_j), \ \ \ j=1,2.
$$
Since $\widehat{T}_1, \widehat{T}_2$, and $L$ commute,
they have common eigenvectors and these vectors
are glued into  a meromorphic function $\Phi(z,\bar{z},P)$ on a
four-sheeted covering
$$
\Gamma({\bf U}) \rightarrow \C: P \in \Gamma({\bf U})
\rightarrow \lambda \in \C,
$$
which ramifies at points where $\widehat{T}_j$ and $L$ are not
diagonalized simultaneously.

To each point $P \in \Gamma({\bf U})$
there corresponds a unique (up to a constant multiple)
Floquet function $\widehat{\varphi}$ with multipliers
$\mu(\gamma_1,P)$ and $\mu(\gamma_2,P)$. If
$\Gamma({\bf U})$ is of finite genus, it is compactified by
four ``infinities'' and $\widehat{\varphi}$ is normalized to
a meromorphic function on $\Gamma({\bf U})$ with exponential
singularities at these ``infinities''.

Notice that if
$\widehat{\varphi} =
(\widehat{\varphi}_1,\widehat{\varphi}_2,\widehat{\varphi}_3,
\widehat{\varphi}_4)^\top$
satisfies (\ref{la-isothermic}) for $\lambda = \mu$ then
\begin{equation}
\sigma(\widehat{\varphi}) =
(\widehat{\varphi}_1,\widehat{\varphi}_2,
-\widehat{\varphi}_3,-\widehat{\varphi}_4)^\top
\label{sigma}
\end{equation}
satisfies (\ref{la-isothermic}) for $\lambda = -\mu$
and this generates an involution
$\sigma:\Gamma({\bf U}) \to\Gamma({\bf U})$,
which descends to an involution of $\C$: $\lambda \to -\lambda$.

By (\ref{sigma}),
the immersion
$$
{\cal M}:\Gamma({\bf U}) \to \C^2 \ \ : \ \ P \to (\mu(\gamma_1,P),
\mu(\gamma_2,P))
$$
descends through the quotient of $\sigma$, i. e.,
${\cal M}:\Gamma({\bf U}) \to \Gamma({\bf U})/\sigma \to \C^2$.

\begin{definition}
The complex curve $\Gamma({\bf U})/\sigma$
is called the spectral curve of an isothermic surface $\Sigma$.
It is said that ${\cal M}(\Gamma({\bf U})/\sigma))$ is the
spectrum of this surface.
\end{definition}

It is again checked by straightforward computations that

\begin{proposition}
If $\widehat{\varphi}$ satisfies (\ref{la-isothermic}), then

1) $\psi = (e^{-\alpha}\varphi_3,e^{-\alpha}\varphi_4)$
satisfies the Dirac equation ${\cal D}\psi = 0$
with the potential
$$
U = \frac{k_1+k_2}{4} e^{\alpha};
$$

2) $\psi^{\ast} = (e^{-\alpha}\varphi_2,e^{-\alpha}\varphi_1)$
satisfies the Dirac equation ${\cal D}\psi^{\ast} = 0$
with the potential
$$
U^{\ast} = \frac{k_2- k_1}{4} e^{\alpha}.
$$
\end{proposition}

We see that the Floquet functions of $(\partial_z - {\bf U})$ and the Dirac
operators ${\cal D}$ with potentials $U$ and $U^{\ast}$ have the same
multipliers and conclude

\begin{theorem}
\label{Tisothermic}
Given an isothermic surface $\Sigma$, the
spectral curve and the spectrum of this isothermic
isothermic surface coincide with

1) a component, of the spectrum of $\Sigma$ as defined in 4.1,
containing both
asymptotic ends where $\mu(\gamma_j) \approx e^{\lambda_+ \gamma_j}$ and
$\mu(\gamma_j) \approx e^{\lambda_- \bar{\gamma}_j}$ as $\lambda_{\pm}
\to \infty$;

2) a component, of the spectrum of the dual surface $\Sigma^{\ast}$
as defined in 4.1,
containing both
asymptotic ends where $\mu(\gamma_j) \approx e^{\lambda_+ \gamma_j}$ and
$\mu(\gamma_j) \approx e^{\lambda_- \bar{\gamma}_j}$ as $\lambda_{\pm}
\to \infty$.

The spectral curve of the isothermic surface
is an irreducible component of the spectral curves of
this surface and its dual as defined in 4.1.
\end{theorem}

Of course, speaking about irreducibility we exclude the case $U \equiv 0$.

Now Pretheorem or, more precisely, its statement about irreducibility of
the spectral curve implies that

{\sl the spectrum and the spectral curve of an isothermic surface
coincide with the spectrum and the spectral curve of this surface and
its dual as defined in 4.1.}

We consider here a general case when surface may be an immersed plane but
not only torus because usually the dual surface to an isothermic torus
is not closed.

In \cite{T2} we introduce a particular case of the
conjecture on the isospectrality of an
isothermic surface and its dual. We conjectured that for surfaces of
revolution, for which the isospectrality is equivalent to coincidence of
all Kruskal--Miura integrals. Theorem \ref{Tisothermic} proves the general
conjecture modulo Pretheorem and, since for surfaces of revolution
Pretheorem  holds (see Example 3 in 3.2 and \cite{T2}), implies the following

\begin{theorem}
A torus of revolution $\Sigma$ and its dual surface $\Sigma^{\ast}$
have the same values of the Kruskal--Miura invariants ${\cal W}_k$.
\end{theorem}

Tori of revolution also explain this passage from $\Gamma({\bf U})$
to $\Gamma({\bf U})/\sigma$.
The spectra of one-dimensional Schr\"odinger
operators related by the Miura transformation (\ref{miura}) coincide
because they are both double covered by the spectrum
of $L$ \cite{EK,T5}.
Moreover the Floquet function $\psi$ of $L$ consists in two components
$\eta_1$ and $\eta_2$ which are the Floquet functions of the Schr\"odinger
operators. In the same manner
the spectra of Dirac operators corresponding to an isothermic
surface and its dual surface are double covered in by the
spectrum of $\widehat{L} =
[\partial_z - \widehat{\bf U}(\lambda)]$ for $\lambda=0$ and
the Floquet function of $\widehat{L}$
consists in the Floquet functions $\psi$ and $\psi^{\ast}$
of the Dirac operators (Proposition \ref{Pisothermic}).

Recall the formulas for the Kruskal--Miura
invariants for surfaces of revolution. Let $\Sigma$ be parameterized by
$x$, a parameter on the rotating curve, and $y$, which is the angle of
rotation and $0 \leq y \leq 2\pi$, and let $z = x+iy$ be a conformal
parameter on $\Sigma$. Let $U$ be the potential of its Weierstrass
representation, which is periodic in $x$ for tori. It also may be defined
for spheres of revolution and
in this case it is fast decaying and defined on the
whole real line \cite{T2}.

The densities of the Kruskal--Miura
integrals for the KdV equation are
$$
R_1 = - q, \ \
R_{n+1} = -R_{nx} - \sum_{k=1}^{n-1} R_k R_{n-k}.
$$
Since $R_{2n}$ are the derivatives of fast decaying functions,
only the integrals
$$
H_n(q) = \int_0^T R_{2n-1}d x
$$
do not vanish identically. Here the integration is taken over $[0,T]$,
where $T$ is the minimal period, for tori and over $\R$ for spheres.
The simplest integrals are
$$
H_1(q) = \int q d x, \ \ H_2(q) = \int q^2 d x, \ \
H_3(q)  = \int (2q^3 - q^2_x) d x.
$$
Put $q = 2 i U_x - U^2$ and define the {\it Kruskal--Miura} invariants as
$$
{\cal K}_l(U) = 2\pi H_l(q).
$$
Notice that ${\cal K}_1$ is the Willmore functional and
others are multiples of ${\cal W}_l$.

\section{Surfaces in the three-sphere}

\subsection{The Dirac equation for surfaces in the
three-sphere}

\medskip

Let $G$ be the Lie group $SU(2)$ and ${\cal G}$ be its
Lie algebra $su(2)$ identified
with the tangent space $T_e G$ to $G$ at the unity $e$.
This Lie algebra is spanned by
\begin{equation}
e_1 =
\left(
\begin{array}{cc}
i & 0 \\ 0 & -i
\end{array}
\right), \ \
e_2 =
\left(
\begin{array}{cc}
0 & 1 \\ -1 & 0
\end{array}
\right), \ \
e_3 =
\left(
\begin{array}{cc}
0 & i \\ i & 0
\end{array}
\right)
\label{basissu}
\end{equation}
which satisfy the commutation relations
$[e_j, e_k] = 2\,\varepsilon_{jkl}\,e_l$.
Take a biinvariant metric on $G$:
$$
\langle \xi, \eta \rangle = - \frac{1}{2}\Tr (\xi\eta), \ \ \
\xi, \eta \in {\cal G} = T_e G,
$$
in which the basis $\{e_1,e_2,e_3\}$ is orthonormal and
$G$ is isometric to the unit $3$-sphere in $\R^4$.

Let $\Sigma$ be a surface immersed into $G$, let $z = x+iy$
be a conformal parameter on $\Sigma$, let
$f:\Sigma \to G$ be the immersion, and let
${\bf I} = e^{2\alpha} dz d\bar{z}$
be the induced metric.

Take the pullback of $TG$ to a  ${\cal G}$-bundle
over $\Sigma$:
${\cal G} \to E = f^{-1}(TG) \stackrel{\pi}{\to} \Sigma$
and the differential
$$
d_{\cal A}:\Omega^1(\Sigma;E) \to \Omega^2(\Sigma;E),
$$
acting on $E$-valued $1$-forms on $\Sigma$ as follows. Let
$$
\omega = f \cdot v dz + f\cdot v^\ast d\bar{z}
$$
where $f\cdot:T_e G \to T_{f(p)}G$
is the left translation by $f(p)$. Then
$$
d_{\cal A} \omega = d'_{\cal A} \omega + d''_{\cal A} \omega
$$
where
$$
d'_{\cal A} \omega = d'_{\cal A} \left(f \cdot v dz \right) = f \cdot
\left(-\bar{\partial}v - \frac{1}{2}[f^{-1}\cdot
\bar{\partial}f,v] \right) dz \wedge d\bar{z},
$$
$$
d''_{\cal A} \omega = d''_{\cal A} \left(f \cdot v^\ast dz \right) = f \cdot
\left(\partial v^\ast + \frac{1}{2}[f^{-1}\cdot
\partial f,v^\ast] \right) dz \wedge d\bar{z}.
$$

By straightforward computations derive that
\begin{equation}
d_{\cal A} (df) = 0.
\label{harmonic1}
\end{equation}

Since $\ast dz = -i dz$ and  $\ast d\bar{z} = i  d\bar{z}$, we have
$$
d_{\cal A} (\ast df) = f\cdot \left( i \bar{\partial}(f^{-1}\cdot\partial f) +
i \partial(f^{-1}\cdot\bar{\partial} f) \right) dz \wedge d\bar{z}.
$$
By the definition of the mean curvature $H$, we have
$$
d_{\cal A} (\ast df) = f\cdot (e^{2\alpha} \tau(f)) dx \wedge dy =
\frac{i}{2} f\cdot(e^{2\alpha} \tau(f)) dz \wedge d\bar{z}
$$
where $\tau(f)$ is the tension vector and
$f\cdot\tau(f) = 2HN$
with $N$ the normal vector:
$f^{-1}\cdot N = -ie^{-2\alpha} [f^{-1}\cdot \partial f,
f^{-1}\cdot\bar{\partial} f]$.
Finally we derive
\begin{equation}
d_{\cal A}(\ast df) = f\cdot \left(H [f^{-1}\cdot\partial f,
f^{-1}\cdot\bar{\partial} f] \right) dz \wedge d\bar{z}.
\label{harmonic2}
\end{equation}
The case $H = 0$ is described by the harmonicity equation
$$
d_{\cal A}(\ast df) = 0.
$$

Put
\begin{equation}
df = f\cdot \left( \Psi dz + \Psi^\ast d\bar{z} \right)
\label{spinor1}
\end{equation}
and rewrite (\ref{harmonic1}) and (\ref{harmonic2}) as
\begin{equation}
\bar{\partial} \Psi - \partial \Psi^\ast + [\Psi^\ast,\Psi] = 0,
\label{harm1}
\end{equation}
\begin{equation}
\bar{\partial} \Psi + \partial \Psi^\ast = iH [\Psi^\ast,\Psi].
\label{harm2}
\end{equation}

Since $f^{-1}\cdot f_x = a^j e_j$ and $f^{-1}\cdot f_y = b^k e_k$
with $a^j, b^k \in \R$, we have
$$
\Psi = \sum_{j=1}^3 Z_j e_j, \ \
\Psi^\ast = \sum_{j=1}^3 \bar{Z}_j e_j
$$
with
$Z_j = (a^j - i b^j)/2$ and $\Psi,\Psi^\ast \in su(2)\otimes \C$.
The induced metric is
$$
e^{2\alpha}dz d\bar{z} = -\frac{1}{2}\Tr[( \Psi dz + \Psi^\ast d\bar{z})^2] =
$$
$$
(Z_1^2 + Z_2^2 + Z_3^2)(dz)^2 + 2(|Z_1|^2 + |Z_2|^2 +|Z_3|^2) dz d\bar{z} +
(\bar{Z}_1^2 + \bar{Z}_2^2 + \bar{Z}_3^2)(d\bar{z})^2
$$
and we conclude that
$$
|Z_1|^2 + |Z_2|^2 +|Z_3|^2 = \frac{1}{2}e^{2\alpha}, \ \ \
Z_1^2 + Z_2^2 + Z_3^2 = 0.
$$
Representing solutions to the latter equation as in 2.2
\begin{equation}
Z_1 = \frac{i}{2}(\bar{\psi}_2^2 + \psi_1^2), \
Z_2 = \frac{1}{2}(\bar{\psi}_2^2 - \psi_1^2), \
Z_3 = \psi_1 \bar{\psi}_2,
\label{spinor2}
\end{equation}
we derive that
$$
e^{\alpha} = |\psi_1|^2 + |\psi_2|^2 .
$$
Now rewriting (\ref{harm1}) in terms of $\psi_j$
and expanding it in the basis $\{e_j\}$, show
that (\ref{harm1}) is equivalent to the system
$$
\bar{\partial}(\psi_1\bar{\psi}_2) - \partial(\bar{\psi}_1 \psi_2) =
i(|\psi_2|^4 - |\psi_1|^4),\ \ \
\bar{\partial}(\psi_1^2) + \partial(\psi_2^2) = 2i\psi_1 \psi_2 e^{\alpha}.
$$
Analogously show that (\ref{harm2}) is equivalent to the system
$$
\bar{\partial}(\psi_1\bar{\psi}_2) +\partial(\bar{\psi}_1 \psi_2) =
H(|\psi_2|^4 - |\psi_1|^4),\ \ \
\bar{\partial}(\psi_1^2) - \partial(\psi_2^2) = 2H\psi_1 \psi_2 e^{\alpha}.
$$
Introduce
$V_1 = -\partial\psi_2 / \psi_1$ and
$V_2 = \bar{\partial} \psi_1 / \psi_2$.
It follows from (\ref{harm1}) that
$\Re V_1 = \Re V_2, \Im V_1 = -e^{\alpha}/2$, and
$\Im V_2 = e^{\alpha}/2$,
and (\ref{harm2}) implies that
$\Re V_1 = \Re V_2 = H e^{\alpha}/2$.

Finally we derive that

\begin{theorem}
For any immersed surface $\Sigma$ is $S^3$ the spinor field
$\psi$ defined by (\ref{spinor1}) and (\ref{spinor2})
satisfies the Dirac equation
$$
{\cal D}^S \psi = 0
$$
with
\begin{equation}
{\cal D}^S =
\left(
\begin{array}{cc}
0 & \partial \\ -\bar{\partial} & 0
\end{array}
\right) +
\left(
\begin{array}{cc}
V & 0 \\ 0 & \bar{V}
\end{array}
\right), \ \ \
V = \frac{1}{2}(H - i)(|\psi_1|^2 + |\psi_2|^2).
\label{dirac-sphere}
\end{equation}
\end{theorem}

This spinor field is unique by its construction and we say that
$\psi$ is the generating spinor for a surface.

Notice that the converse is not always true as in the case of Theorems 1 and 2.
Indeed, not to any solution of ${\cal D}^S\psi = 0$
corresponds a surface: a solution related to a surface has to satisfy
an additional condition
$$
|\psi_1|^2 + |\psi_2|^2 = - 2 \Im V.
$$

It is easy to check that if ${\cal D}^S\psi = 0$, then
${\cal D}^S\varphi = 0$ with $\varphi = (\bar{\psi}_2,-\bar{\psi}_1)^{\top}$.

Let us write the complete system of the Gauss--Weingarten equations.
Recall that the Hopf differential equals
$A dz^2 = \langle f_{zz}, N \rangle\, dz^2$
and, since the metric is left-invariant, we have
$A = \langle f^{-1}f_{zz}, f^{-1} N \rangle$.
Now $\Psi = f^{-1}f_z$ and
$$
f^{-1}f_{zz} = \Psi_z + \Psi^2,
$$
where
\begin{equation}
\Psi =
\left(
\begin{array}{cc}
i Z_1 & Z_2 + i Z_3 \\
-Z_2 + i Z_3 & -i Z_1
\end{array}
\right), \ \ \
\Psi^{\ast} =
\left(
\begin{array}{cc}
i \bar{Z}_1 & \bar{Z}_2 + i \bar{Z}_3 \\
-\bar{Z}_2 + i \bar{Z}_3 & -i \bar{Z}_1
\end{array}
\right).
\label{Psi}
\end{equation}
We have $\Psi^2 = (Z_1^2 + Z_2 ^2 + Z_3^2) e_1$ and, since $z$ is a conformal
parameter, $\Psi^2 = 0$.
Therefore as in 2.2 the Hopf differential takes the same form
$$
A\, d z^2 = (\psi_{1z}
\bar{\psi}_2 - \bar{\psi}_{2z}\psi_1) d z^2.
$$
We also have
$$
\alpha_z e^{\alpha} = \bar{\psi}_1\psi_{1z} + \psi_2 \bar{\psi}_{2z}
$$
and finally write down the Gauss--Weingarten equations for an immersed surface
in $S^3$ as
\begin{equation}
\left[\frac{\partial}{\partial z} -
\left(\begin{array}{cc}
 \alpha_z & A e^{-\alpha} \\
-V & 0
\end{array}
\right)\right]\psi =
\left[\frac{\partial}{\partial \bar{z}} -
\left(\begin{array}{cc}
0 & \bar{V} \\
-\bar{A}e^{-\alpha} & \alpha_{\bar{z}}
\end{array}
\right)\right]\psi = 0.
\label{w-gw-sphere}
\end{equation}
The compatibility conditions are the Codazzi equations
\begin{equation}
\alpha_{z\bar{z}} + |V|^2 - |A|^2e^{-2\alpha} = 0, \ \ \
A_{\bar{z}} = (\bar{V}_z - \alpha_z\bar{V})e^{\alpha}.
\label{codazzi-sphere}
\end{equation}

{\sl Examples.}

{\sl 1) The Clifford torus.}
This torus in $\R^4$ is defined by the equations
$$
(x^1)^2 + (x^2)^2 = (x^3)^2 + (x^4)^2 = \frac{1}{2}
$$
where $(x^1,\dots,x^4) \in \R^4$. It is immersed
into $SU(2)$ by the formula
$$
f(x,y) = \frac{1}{\sqrt{2}}
\left(
\begin{array}{cc}
e^{ix} & e^{iy} \\ -e^{-iy} & e^{-ix}
\end{array}
\right)
$$
and has a conformal type of the square torus, $(x,y) \in
\Z^2/2\pi\Z^2$.  Compute that
$$
\Psi = \frac{1}{4}((1+i)e_1 +
(1-i)\sin (x-y) e_2 + (1-i) \cos (x-y) e_3),
$$
$$
\psi_1 =
\sqrt{\frac{1-i}{2}} \sin\left(\frac{x-y}{2} - \frac{\pi}{4}\right), \ \ \
\psi_2 =
\sqrt{\frac{1+i}{2}} \cos\left(\frac{x-y}{2} - \frac{\pi}{4}\right),
$$
$$
e^{\alpha} = \frac{1}{\sqrt{2}}, \ \ \  V = -\frac{i}{2\sqrt{2}}, \ \ \
A = \frac{1}{4}.
$$

{\sl 2) Minimal tori in $S^3$.}
In this case we have $V = -ie^{\alpha}/2$ and derive
from (\ref{codazzi-sphere})
that $A_{\bar{z}} = 0$. This means that the Hopf differential is holomorphic
and as in the case of CMC tori in $\R^3$ we conclude that it is constant
and by rescaling conformal parameter achieve $A = 1/2$.
The case $A=0$ is also excluded for tori: it is realized by the equatorial
$S^2$-spheres in $S^3$ (complete CMC surfaces in $\R^3$ with $A=0$ are
the round spheres). The first
equation from (\ref{codazzi-sphere}) is
$$
u_{z\bar{z}} + \sinh u = 0, \ \ \ \ u=2\alpha.
$$

For CMC tori in $S^3$ the Hopf differential is also holomorphic and they are
described in the same manner.

Now it is clear that the analogs of Theorems 1, 2, and 3 hold for surfaces
in $S^3$. We again have the same spinor bundles over constant curvature
surfaces.

\begin{definition}
Given a torus, represented in $S^3$ via $\psi$ satisfying (\ref{w-gw-sphere}),
the spectral curve $\Gamma^S$ of the operator ${\cal D}^S$ with the potential
(\ref{dirac-sphere}) is called the spectral curve of the torus. Given
in addition
a basis $\gamma_1$,$\gamma_2$ for $\Lambda$, the period lattice for $U$,
the image of the multiplier map
$$
{\cal M}:Q_0({\cal D}^S)/\Lambda^{\ast} \to \C^2 \ \ \ :
\ \ \
{\cal M}(k) = (e^{2\pi i \langle k,\gamma_1 \rangle},
e^{2\pi i \langle k,\gamma_2 \rangle})
$$
is called the spectrum of the torus in $S^3$.
\end{definition}

{\sl Example.} Let $\Sigma$ be the Clifford torus. Then
$V = -i/2\sqrt{2}$ and the Floquet eigenfunctions
$\psi=(\psi_1,\psi_2)^{\top}$ satisfy the equation
$$
\left(\partial \bar{\partial} + \frac{1}{8}\right)\psi_j = 0, \ \ j=1,2.
$$
We derive that the general Floquet function is
$$
\psi(z,\bar{z},\lambda) =
\left(e^{\lambda z - \frac{1}{8\lambda}\bar{z}},
\frac{i}{2\sqrt{2}\lambda}e^{\lambda z -
\frac{1}{8\lambda}\bar{z}}\right)^{\top}
$$
and find the spectrum as the image of the multiplier map
$$
\lambda \in \C \setminus \{0\} \to
\left(e^{2\pi(\lambda - \frac{1}{8\lambda})},
e^{2\pi i(\lambda + \frac{1}{8\lambda})} \right).
$$
The spectral curve is the two-sphere, i.e., the punctured
plane $\C \setminus \{0\}$ compactified by the points $\lambda = 0,\infty$.
This implies that the spectral genus of the Clifford torus equals zero.

We shall not discuss the spectra of tori in $S^3$ in detail but
only mention that Pretheorem also has to hold for them.

\subsection{The Hitchin system}

Let us compare the previous computations with the Hitchin theory of
harmonic tori in the $3$-sphere \cite{Hitchin}.
For Riemannian manifolds $N$ and $M$ a mapping $f:N \to M$ is called
{\it harmonic} if it satisfies the equations
$$
d_{\cal A}(df) = d_{\cal A}(\ast df) =0
$$
where ${\cal A}$ is the pullback of the Levi-Civita connection on $TM$ and
the Hodge operator $\ast$ is taken with respect to the metric on
$N$. If $f$ is an immersion and the metric on
$N$ is the induced metric, then $f(N)$ is a minimal submanifold.

Let $N$ be an immersed surface $\Sigma$ with the induced metric and
$M=G$ be Lie group with a biinvariant metric.
We adopt the notation from 6.1.

The harmonicity equations take the form
$$
d_{\cal A} (\Psi dz + \Psi^{\ast} d\bar{z}) = 0,
$$
$$
d_{\cal A}(\ast (\Psi dz + \Psi^{\ast} d\bar{z})) = -i
d_{\cal A} (\Psi dz - \Psi^{\ast} d\bar{z}) = 0.
$$
They describe minimal surfaces in $S^3$ and are rewritten as
\begin{equation}
\bar{\partial} \Psi - \partial \Psi^\ast + [\Psi^\ast,\Psi] =
\bar{\partial} \Psi + \partial \Psi^\ast = 0.
\label{hitchin1}
\end{equation}
Following \cite{Hitchin} put
$$
\Phi = \frac{1}{2}\Psi,  \ \ \ \Phi^{\ast} =  - \frac{1}{2}\Psi^{\ast}
$$
and rewrite (\ref{hitchin1}) as the {\it Hitchin system}
\begin{equation}
d^{\prime\prime}_{\cal A} \Phi = 0, \ \ \ F_{\cal A} = d_{\cal A}^2 =
[\Phi,\Phi^{\ast}] = 0,
\label{hitchin2}
\end{equation}
where $F_A$ is the curvature of the connection
$$
d_{\cal A}:\Omega^p(\Sigma;f^{-1}TG) \to \Omega^{p+1}(\Sigma;f^{-1}TG)
$$
and the formula means that $d_{\cal A}^2$ coincides with the multiplication
by $F_{\cal A}$.
The system (\ref{hitchin2}) describes general harmonic mappings of
surfaces in $S^3$ (when the metric on the surface is not necessarily
induced) in terms of a connection $\cal A$ associated to the harmonic map
and the {\it Higgs field} $\Phi$.

The equation $d_{\cal A}df = 0$ is equivalent to
$$
\bar{\partial} \Psi - \partial \Psi^\ast + [\Psi^\ast,\Psi] = 0
$$
and means that the connection
${\cal A}= (\partial + \Psi, \bar{\partial} + \Psi^{\ast})$
on $f^{-1}TG$ is flat, which is evident from its construction.
However the second of the equations (\ref{hitchin1})
implies that this connection is extended to an analytic family
of flat connections
$$
{\cal A}_{\lambda} = \left(\partial + \frac{1+\lambda^{-1}}{2}\Psi,
\bar{\partial} + \frac{1+\lambda}{2}\Psi^{\ast}
\right)
$$
where ${\cal A} = {\cal A}_1$ and $\lambda \in \C\setminus \{0\}$.
This commutation representation with a spectral parameter
was found by Pohlmeyer \cite{Pohl}
for harmonic maps into $SU(2)$ and
later developed by Mikhailov and Zakharov for the case, when the target space
is not a group but a symmetric space $S^2$ \cite{ZM}.
Finally these two papers gave rise to the ``integrability'' part
of the modern theory of harmonic maps \cite{U,Hitchin,BFPP,Guest}.

For harmonic tori Hitchin introduced spectral curves and showed
that they are of finite genus \cite{Hitchin}. Their construction
is as follows.

Let $\Sigma$ be a harmonic torus in $S^3$.
For any $\lambda \in \C \setminus \{0\}$ we have a flat $Sl(2,\C)$ connection.
Fix a basis $\{\gamma_1,\gamma_2\}$ for $H_1(\Sigma)$. For $\gamma_1$ and
$\gamma_2$
define matrices $H(\lambda), \widetilde{H}(\lambda) \in SL(2,\C)$
which describe the monodromies of ${\cal A}_{\lambda}$ along closed loops
realizing $\gamma_1$ and $\gamma_2$.
These matrices commute and have joint eigenfunctions
$\varphi(\lambda,\mu)$ where $\mu$ is a root of
the characteristic equation
for $H(\lambda)$
$$
\mu^2 - \Tr H(\lambda) + 1 = 0
$$
and therefore there is a Riemann surface
on which the eigenvalues
$$
\mu_{1,2} = \frac{1}{2}\left(
\Tr H(\lambda) \pm \sqrt{\Tr^2 H(\lambda) -4}\right)
$$
are defined.
The complex curve $\Gamma$, which is a two-sheeted covering of $\C P^1$,
ramifying at the odd zeros of the function $(\Tr^2 H(\lambda) -4)$
and at $0$ and $\infty$, is called the {\it spectral curve of a
harmonic torus} in $S^3$.

On $\Gamma$ the eigenvalues of
$H(\lambda)$ paste into a single-valued function $\mu$ with singularities
at $0$ and $\infty$.
Moreover the joint eigenfunctions of $H(\lambda)$ and $\widetilde{H}(\lambda)$
paste into a vector function $\varphi$
meromorphic on $\Gamma \setminus \{0,\infty\}$.

We shall show that in the case, when the harmonic tori is
an immersed tori in $S^3$ with the induced metric, i. e., in the situation of
6.1, the spectral curve of Hitchin is the same as the spectrum of the torus
as defined in 6.1.

Let $f:\Sigma \to S^3$ be an immersion of a minimal torus
and $\Psi = f^{-1}f_z, \Psi^{\ast}= f^{-1}f_{\bar{z}}$.
Let the surface be defined by a spinor $\psi$.

The Hitchin eigenfunction $\varphi(\lambda,\mu)$ satisfies the equations
$$
\left[ \partial + \frac{1+\lambda}{2}\Psi \right] \varphi =
\left[ \bar{\partial} +
\frac{1+\lambda^{-1}}{2}\Psi^{\ast} \right] \varphi = 0.
$$
Take the matrix
\begin{equation}
L =
\left(
\begin{array}{cc}
\bar{a} & -\bar{b} \\ b & a
\end{array}
\right)
\label{compat}
\end{equation}
with
$a = (-i\bar{\psi}_1 + \psi_2)/\sqrt{2},
b = (- i \psi_1 + \bar{\psi}_2)/\sqrt{2}$.
By (\ref{spinor2}) and (\ref{Psi}), compute that
$$
L^{-1} \Psi L =
e^{\alpha}
\left(
\begin{array}{cc}
0 & 1 \\ 0 & 0
\end{array}
\right), \ \ \ \ \
L^{-1} \Psi^{\ast} L =
e^{\alpha}
\left(
\begin{array}{cc}
0 & 0 \\ -1 & 0
\end{array}
\right).
$$
We also have
$$
L^{-1}L_z =
\left(
\begin{array}{cc}
\alpha_z & -iV \\ -iAe^{-\alpha} & 0
\end{array}
\right), \ \ \ \ \
L^{-1}L_{\bar{z}} =
\left(
\begin{array}{cc}
0 & -i\bar{A}e^{-\alpha} \\ -i\bar{V} & \alpha_{\bar{z}}
\end{array}
\right).
$$
The vector function $L^{-1} \varphi$ satisfies the equations
$$
\left[
\partial +
\left(
\begin{array}{cc}
\alpha_z & -iV \\ -iAe^{-\alpha} & 0
\end{array}
\right) + \frac{1+\lambda}{2}e^{\alpha}
\left(
\begin{array}{cc}
0 & 1 \\ 0 & 0
\end{array}
\right)
\right] L^{-1}\varphi = 0,
$$
$$
\left[
\bar{\partial} +
\left(
\begin{array}{cc}
0 & -i\bar{A}e^{-\alpha} \\ -i\bar{V} & \alpha_{\bar{z}}
\end{array}
\right) +
\frac{1+\lambda^{-1}}{2}
e^{\alpha}
\left(
\begin{array}{cc}
0 & 0 \\ -1 & 0
\end{array}
\right)
\right] L^{-1}\varphi = 0.
$$
These two equations are compatible only for minimal tori, which
are described by the condition
$$
V = - \frac{ie^{\alpha}}{2}.
$$
For $\widetilde{\varphi} = e^{\alpha}L^{-1}\varphi$ we derive that
$$
\partial \widetilde{\varphi}_1 + \frac{\lambda}{2}e^{\alpha}
\widetilde{\varphi}_2 = 0,
\ \ \ \
\bar{\partial}\widetilde{\varphi}_2 -
\frac{1}{2\lambda}e^{\alpha}\widetilde{\varphi}_1 = 0.
$$
Put $\widetilde{\psi}_1 = i \lambda \widetilde{\varphi}_2,
\widetilde{\psi}_2 = \widetilde{\varphi}_2$
and notice that $\widetilde{\psi}$ satisfies
(\ref{dirac-sphere}):
$$
\left[
\left(
\begin{array}{cc}
0 & \partial \\
-\bar{\partial} & 0
\end{array}
\right)
+
\left(
\begin{array}{cc}
V & 0 \\
0 & \bar{V}
\end{array}
\right)
\right]
\widetilde{\psi} = 0
\ \ \ \
\mbox{with $V = -ie^{\alpha}/2$}.
$$
As in the proofs of Theorems 5 and 6
we conclude

\begin{theorem}
Given a minimal torus $f:\Sigma \to S^3$,
the Hitchin eigenfunction $\varphi(\lambda,\mu)$ by the transformation
$$
\left(
\begin{array}{c}
\varphi_1 \\
\varphi_2
\end{array}
\right) \to
\left(
\begin{array}{c}
\widetilde{\psi}_1 \\
\widetilde{\psi}_2
\end{array}
\right)
=
e^{\alpha}
\left(
\begin{array}{cc}
0 & i \lambda \\
1 & 0
\end{array}
\right) \cdot
L^{-1}
\cdot
\left(
\begin{array}{c}
\varphi_1 \\
\varphi_2
\end{array}
\right)
$$
is mapped to a Floquet function $\widetilde{\psi}$
of ${\cal D}^S$.
There is the mapping of
the eigenvalues $\mu$ and $\widetilde{\mu}$
of $\varphi$
with respect to the monodromy operators $H(\lambda)$ and
$\widetilde{H}(\lambda)$ to the multipliers of $\widetilde{\psi}$:
\begin{equation}
(\mu, \widetilde{\mu}) \to
((-1)^{\varepsilon(\gamma_1)}\mu,
(-1)^{\varepsilon(\gamma_2)}\widetilde{\mu})
\label{multis}
\end{equation}
where $(-1)^{\varepsilon(\gamma_1)},(-1)^{\varepsilon(\gamma_2)}$
are the multipliers of the spinor $\psi$ generating a minimal torus.

This mapping (\ref{multis}) establishes a biholomorphic equivalence between
the Hit\-chin spectral curve of a minimal torus and the connected component
of the spectrum of this torus as defined in 6.1. This connected
component contains both asymptotic ends near which $\widetilde{\psi} \approx
(e^{\lambda_+}z,0)^{\top}$ or $\widetilde{\psi} \approx
(0,e^{\lambda_-}\bar{z})^{\top}$.
\end{theorem}

If Pretheorem holds for ${\cal D}^S$, then the spectrum is irreducible and
therefore (\ref{multis}) establishes a biholomorphic equivalence of both
spectra.

\section{Conformal invariance of the spectra of tori}

\subsection{The M\"obius group}

We consider $\R^4$ as the set of matrices
\begin{equation}
a = \left(
\begin{array}{cc}
x^4 + i x^1 & x^2 + i x^3 \\
-x^2 + i x^3 & x^4 - ix^1
\end{array}
\right), \ \ \ \ \ x^1,x^2,x^3,x^4 \in \R,
\label{matrixrep}
\end{equation}
and consider $\R^3$ as a subset described by $x^4=0$.
The unit sphere $S^3=SU(2)$ is defined by the equation
$$
|x| = 1.
$$
Take the north pole $P = (0,0,0,1)$ and denote by $\pi$ the stereographic
projection of $S^3$ to $\R^3 = \{x^4=0\}$ from $P$:
$$
\pi: a \to
\frac{1}{1-x^4}
\left(
\begin{array}{cc}
i x^1 & x^2 + i x^3 \\
-x^2 + i x^3 & -ix^1
\end{array}
\right) = (1+a)(1-a)^{-1}.
$$
The inverse mapping is
$$
\pi^{-1}: b \to (b-1)(b+1)^{-1}.
$$
This mapping $\pi$ establishes a conformal equivalence between $S^3$ and
$\R^3$ compactified by a point at infinity, i.e., by $\pi(P) = \infty$.

The group of conformal transformations of $\R^3 \cup \infty$ is isomorphic to
$O^+(4,1)$, the subgroup of $O(4,1)$ formed by isochronic transformations.
The geometric picture is as follows.
Let $\R^{1,4}$ be a $5$-dimensional pseudo-Euclidean space with the
metric
$$
\langle x,y \rangle_{1,4} = x^0 y^0 - \sum_{j=1}^4 x^j y^j.
$$
The $4$-dimensional hyperbolic space ${\cal H}^4$
is embedded into $\R^{1,4}$
as the upper half of a hyperboloid:
$\langle x, x \rangle_{1,4} = 1, x^0 > 0$,
with the metric on tangent vectors
$\langle \xi, \xi \rangle = - \langle \xi, \xi \rangle_{1,4}$.
The group of isometries of ${\cal H}^4$ is $O^+(4,1)$ and it
acts on $S^3$, the absolute of ${\cal H}^4$, by conformal transformations.

By the Liouville theorem,
the group $O^+(4,1)$ of conformal transformations is generated by

1) isometries of $\R^3$;

2) inversions with centers in $x_0 \in \R^3$: $x \to \frac{x-x_0}{|x-x_0|^2}$;

3) homotheties: $x \to \lambda x$, $\lambda \in \R \setminus \{0\}$.

Any conformal transformation of $\bar{\R}^3 =
\R^3 \cup \{\infty\}$ which
preserves $\infty$ is a composition of isometries and homotheties.
Notice that in terms of (\ref{matrixrep}) the inversion of $\R^3$ centered at
$x_0$ looks simply as
$x \to (x_0-x)^{-1}$.

We think that

{\sl the spectrum of a torus in $S^3$ which is stereographically
projected into a torus in $\R^3$  and the spectrum of this
projection coincide.}

We can not prove it now but would like to notice that this statement
easily implies conformal invariance of both spectra: it is clear that
the spectrum of a torus in $\R^3$ is invariant under translations
of the torus and the spectrum of a torus in $S^3$ is invariant under
rotations in $S^3$. However the stereographic projection converts
rotations in $S^3$ into conformal transformations of $\R^3$ which together
with translations and homotheties generate
the conformal group $O^+(4,1)$.
The same holds with translations of $\R^3$ whose compositions with
the projection generate together with rotations
the group of conformal transformations of $S^3$.

\subsection{Conformal invariance of the spectra for
isothermic tori in $\R^3$}

\begin{theorem}
Let $\Sigma$ be an isothermic torus in $\R^3$ and let $F:
\bar{\R}^3 \to \bar{\R}^3 $ be a conformal transformation which maps
$\Sigma$ into a torus $F(\Sigma)$ lying in $\R^3$.
Then the spectrum of an isothermic torus $\Sigma$ coincides with
the spectrum of $F(\Sigma)$.
\end{theorem}

Notice that the spectrum of an isothermic torus is defined as a component of
the general spectrum which contains
the asymptotic ends where $\mu(\gamma_j) \approx e^{\lambda_+ \gamma_j}$
and $\mu(\gamma_j) \approx e^{\lambda_j \bar{\gamma}_j}$ as
$\lambda_{\pm} \to \infty$.
Pretheorem states that the general spectrum is irreducible and therefore
coincides with this component.

{\sl Proof.}
By Theorem 6, the spectrum of an isothermic torus and its dual
isothermic surface coincide. The potential of the dual surface equals
$$
U^{\ast} = \frac{k_2-k_1}{4}e^{\alpha}
$$
and, by the Blaschke theorem, the density of the
Willmore functional
$$
\left(\frac{k_2-k_1}{2}\right)^2 d\mu =
4 \left(U^{\ast}\right)^2 d x \wedge d y
$$
is invariant under conformal transformations of $\bar{\R}$.
As known, conformal transformations maps isothermic surfaces into
isothermic ones.

Let $z$ be a conformal parameter on $\Sigma$ which is mapped into a
conformal parameter on $F(\Sigma)$ and $V$ be the potential of $F(\Sigma)$
with respect to this parameter. We see that, by the Blaschke theorem,
 $V^2 = (U^{\ast})^2$ and therefore $V = \pm U^{\ast}$.

It is clear that the spectra of the Dirac operators whose potentials differs
by sign coincide. Now we derive that the spectra of the isothermic tori
$\Sigma$ and $F(\Sigma)$ coincide with the spectrum of the isothermic surface
with the potential $U^{\ast}$. This proves the theorem.

\end{document}